\documentclass[a4paper,12pt]{article}

\usepackage{amsmath,amsthm,amssymb,enumerate}
\usepackage{xcolor}
\usepackage{titling}
\usepackage{float}
\usepackage{multirow}
\usepackage{longtable}
\usepackage{diagbox}
\usepackage{soul}
\usepackage{cancel}
\usepackage{mathrsfs}
\usepackage[affil-it]{authblk}
\usepackage{gensymb}
\usepackage{enumitem}
\usepackage[square,sort&compress,comma,numbers]{natbib}
\usepackage{bm}
\usepackage{hyperref}
\usepackage{subfig}
\usepackage{graphicx}
\usepackage[margin=1.7cm]{geometry}
\usepackage{tikz}
\usepackage{esint}
\usepackage{caption}
\usetikzlibrary{shapes,calc}
\usepackage{verbatim}
\usepackage{array}
\usepackage{bm}

\usepackage{newtxtext,newtxmath}

\usepackage{multirow}
\usepackage{marginnote}
\chardef\bslash=`\\ %

\newtheorem{defn}{Definition}
\theoremstyle{definition}
\theoremstyle{rem}

\newtheorem{rem}{Remark}[section]
\newtheorem{example}{Example}[section]

\numberwithin{equation}{section}

\newcommand{\bT}{\mathbb T}

\newcommand{\cE}{\mathcal E}

\newcommand{\cT}{\mathcal{T}}

\newcommand{\cM}{{\rm M}}

\newcommand{\vket}{von K\'{a}rm\'{a}n equations}

\newcommand{\fl}{\quad \text{for all}\:}

\newcommand{\half}{\frac{1}{2}}

\newcommand{\dx}{{\rm\,dx}}
\newcommand{\ds}{\!{\rm ds}}

\newcommand{\dg}{{\rm dG}}

{\algorithm}%
{\endalgorithm}%

\theoremstyle{definition}

\numberwithin{equation}{section}

\newcommand{\bV}{\text{\bf V}}

\newcommand{\bv}{\boldsymbol{v}}

\newcommand{\Stb}{{\boldsymbol{P}_2}(\cT)}

\newcommand{\ip}{{\rm IP}}

\newcommand{\trinl}{\ensuremath{|\!|\!|}}
\newcommand{\trinr}{\ensuremath{|\!|\!|}}

\newcommand{\TT}[1]{\mathbb T\left(#1\right)}

\DeclareMathOperator{\E}{\mathcal{E}}
\newcommand{\M}{\mathrm{M}}

\newcommand{\T}{\mathcal{T}}

\newcommand{\trb}[1]{|\!|\!|#1|\!|\!|}

\newcommand{\id}{\mathrm{id}}

\newcommand{\LP}{\Lambda_{\rm P}}

\renewcommand{\div}{\mathrm{div}}

\newcommand{\ua}{{u_{1}}}
\newcommand{\ub}{{u_{2}}}

\newcommand{\vha}{{v_{h,1}}}
\newcommand{\vhb}{{v_{h,2}}}
\newcommand{\vhl}{{v_{h,\ell}}}

\theoremstyle{remark}
\theoremstyle{plain}
\newtheorem{theorem}{Theorem}[section]
\newtheorem{lemma}[theorem]{Lemma}
\usepackage{totcount}
\newtotcounter{nconst}
\makeatletter
\newcount\rc@count
\let\rc@clearconstantlist\empty

\newcommand\rc@clearconstant[1]{\global\expandafter\let\csname rc@const@#1\endcsname\undefined}

\newcommand\resetconstants[1]{%
    \def\rc@constname{#1}% Set the new base name of the constants to the argument
    \global\rc@count=1\relax % Reset the constant counter to 1
    \bgroup 
        \let\\\rc@clearconstant % map over the list of constants that have been defined, clearing each of them.
        \rc@clearconstantlist
        \global\let\rc@clearconstantlist\empty % Globally empty the list of constants.
    \egroup
}
\newcommand\const[1]{%
    \@ifundefined{rc@const@#1}{%
        \expandafter\xdef\csname rc@const@#1\endcsname{%
           \noexpand\rc@useconst{\rc@constname}{\the\rc@count}%
        }%
        \g@addto@macro\rc@clearconstantlist{\\{#1}}%
        \global\advance\rc@count1\relax
    }{}%
	\setcounter{nconst}{\the\rc@count-1}
    \csname rc@const@#1\endcsname
}

\newcommand\rc@useconst[2]{\ensuremath{#1_{#2}}}
\makeatother

	\resetconstants{C}

\newcommand{\arcangle}{%
  \mathord{<\mspace{-7mu}\mathrel{)}\mspace{3mu}}%
}

\newcommand{\dislin}{a_h}
\newcommand{\rop}{R}
\newcommand{\sop}{S}
\newcommand{\lamr}{\Lambda_{\rm R}}
\newcommand{\lams}{\Lambda_{\rm S}}

\def \R{{{\Bbb R}}}

\allowdisplaybreaks

\def\R{\mathbb{R}}

\def\O{\Omega}

\def\bv{{\mathbf V}}

\def\pw{{\rm {pw}}}

\def\cE{{\mathcal{E}}}

\def\jump#1{\left[ #1 \right]}

\newcommand{\be}{\begin{equation}}
\newcommand{\ee}{\end{equation}}

\newcommand{\pl}{\delta_\zeta}
\usepackage{mathtools}  
\usepackage[normalem]{ulem}
\definecolor{violet}{rgb}{0.580,0.,0.827}

\usepackage[normalem]{ulem}
\newcounter{corr}
\definecolor{violet}{rgb}{0.580,0.,0.827}
\newcommand{\corr}[3]{\typeout{Warning : a correction remains in page
		\thepage}
	\stepcounter{corr}        
	{\color{red}\ifmmode\text{\,{\ensuremath{#1}}\,}\else{#1}\fi}
	{\color{blue}#2}
	{\color{violet} #3}}

\newcounter{changeto}
\newcommand{\changeto}[2]{\typeout{Warning : a correction remains in page
		\thepage}
	\stepcounter{changeto}        
	{\color{blue}\ifmmode\text{\,\sout{\ensuremath{#1}}\,}\else\sout{#1}\fi}
	{\color{red}#2}}

 \usepackage{xparse}

   \newcounter{const}
\NewDocumentCommand{\constant}{o}
 {%
  \IfValueTF{#1}%
  {C_{#1}}%
  {\refstepcounter{const}%
  C_{\theconst}}%
 }

\title{A posteriori error control {for} fourth-order semilinear problems with quadratic nonlinearity }
\author{Carsten Carstensen\footnote{Department of Mathematics, 
Humboldt-Universit\"{a}t zu Berlin, 10099 Berlin, Germany.
Distinguished Visiting Professor,  Department of Mathematics, Indian Institute of 
Technology Bombay, Powai, Mumbai-400076, India.  cc@math.hu-berlin.de}        
\quad\text{and}\quad Benedikt Gr\"a{\ss}le\footnote{Department of Mathematics, 
Humboldt-Universit\"{a}t zu Berlin, 10099 Berlin, Germany.
graesslb@math.hu-berlin.de}\quad\text{and}\quad
Neela Nataraj\footnote{Department of Mathematics,~Indian Institute of Technology Bombay,~Powai, Mumbai 400076,~India.~neela@math.iitb.ac.in}
}

\begin{document}

	\maketitle
\abstract{
A general  a posteriori error analysis applies to  five lowest-order finite element methods for  two fourth-order semi-linear problems with trilinear non-linearity and a general source.   A quasi-optimal smoother extends the source term to the discrete trial space, and more importantly, modifies the trilinear term in the stream-function vorticity formulation of the incompressible 2D Navier-Stokes and the  \vket.  This enables the first efficient and reliable a posteriori error estimates for the 2D Navier-Stokes equations in the stream-function vorticity formulation for
Morley, two discontinuous Galerkin, $C^0$ interior penalty, and WOPSIP discretizations with piecewise quadratic polynomials.}
	


	\section{Introduction}
	This paper discusses an abstract a posteriori error analysis for fourth-order semilinear problems and
	its applications to the incompressible 2D Navier-Stokes equations and the von K\'arm\'an equations.
The continuous problem in this paper seeks a regular root $u$ 
in a Banach space $X$   to $N \in C^1(X;Y^*)$ for
\begin{equation}\label{eqn:cont}
N (x):= a(x,\bullet) + {\Gamma}(x,x,\bullet)- F \text{ for } x \in X. 
\end{equation}
The bilinear form $a(\bullet,\bullet)$ in \eqref{eqn:cont} corresponds to a weak form of the biharmonic operator, the trilinear form $\Gamma(\bullet,\bullet,\bullet)$ represents a quadratic nonlinearity, and $F$ is the general source term in $Y^*$; for instance, for the 2D Navier-Stokes equations in the stream-function vorticity formulation and for the von K\'{a}rm\'{a}n plates. 
The nonconforming discretization of \eqref{eqn:cont} with a piecewise application of the differential operators in the weak forms for Morley finite element functions \cite{CCGMNN_semilinear} allows for a priori convergence results. But their a posteriori error analysis so far was not satisfactory for the stream-function vorticity formulation of the incompressible 2D Navier-Stokes equations \cite{CCGMNN_semilinear, kim_morley_2021}: the efficiency analysis is excluded in \cite{CCGMNN_semilinear} and merely partial in  \cite[Remark
 4.11]{kim_morley_2021}. 
 For the (generalised) Morley interpolation operator $I_\M$, the companion operator $J$ \cite{aCCP, ccnnlower2022,CCBGNN22}, the smoother $J_h=JI_\M=S=Q$,  the  choice $R \in \{{\rm id}, I_{\M},  JI_\M\}$, 
the discrete problem \cite{CCNNGRDS22} seeks an approximation $u_h$ to a regular root $u$ to \eqref{eqn:cont}   in  a finite-dimensional space $X_h$ as a root of 
\begin{equation}\label{eqn:discrete}
N_h(u_h) := \dislin(u_h,\bullet) + {\Gamma}_{\pw}(\rop u_h,\rop u_h,\sop\bullet)- F_h.%
\end{equation}
 The bilinear form $a_h(\bullet,\bullet)$  discretizes $a(\bullet,\bullet)$, for instance, with Morley \cite{DG_Morley_Eigen,CCDG14,CCDGJH14},  discontinuous Galerkin (dG) \cite{Georgoulis2009, Georgoulis2011,CCGMNN18}, $C^0 $ interior penalty (IP) \cite{BS05,BS_C0IP_VKE},  and WOPSIP \cite{BrenGudiSung10} schemes;  the trilinear form $\Gamma_{\pw}(\bullet,\bullet,\bullet)$ discretizes $\Gamma(\bullet,\bullet,\bullet)$ by the piecewise action of the differential operators,  $F_h= F \circ Q$ approximates $F$,  and $R$, $S$ denote quasi-optimal smoothers in the spirit of \cite{veeser1,veeser2,veeser3,CCBGNN22,ccnnlower2022,CCNNGRDS22}. The innovative point  in \eqref{eqn:discrete} is the application of smoothers $R$ and $S$ in the nonlinearity $\Gamma_{\pw}(Ru_h,R u_h, Sy_h)$.
The prequel \cite{CCNNGRDS22} establishes an a priori analysis of this class of lowest-order finite element methods and a source term $F \in Y^*$ with the {first} best-approximation result for $S=Q=J_h$, namely 
\begin{align}\label{eqn:QO}
    \|u - u_h \|_{\widehat{X}} \le 
C_{\rm qo} \min_{x_h \in X_h} \|u - x_{h}\|_{\widehat{X}}.\tag{\bf QO}
\end{align}
Here and throughout this paper, the Banach spaces $(X_h,\|\bullet\|_{X_h})$ and $(X,\|\bullet\|_X$) are contained in a common superspace $(\widehat X, \|\bullet\|_{\widehat X})  $ with a norm $\|\bullet\|_{\widehat X}$ that extends $\|\bullet\|_{X_h}=(\|\bullet\|_{\widehat X})|_{X_h}$ and $\|\bullet\|_{X}=(\|\bullet\|_{\widehat X})|_{X}$.
  This paper presents the first reliable and efficient  a posteriori error analysis for those schemes and includes the {\it first} reliable and efficient a posteriori estimates for the lowest-order finite element schemes for the  2D Navier-Stokes equation in the stream-function vorticity formulation.

\smallskip \noindent

\smallskip \noindent 
Section~\ref{sec:Abstract a posteriori error control} 
introduces an abstract framework of an a posteriori error control in Banach spaces $X$ and $Y$ as in \eqref{eqn:cont}-\eqref{eqn:discrete} that applies below  to five second-order schemes, namely the Morley, two dG, $C^0$IP, and WOPSIP. The outcome allows for rough source terms $F \in Y^*$ and provides reliable and efficient error control by the sum of three contributions. Given an approximation $v_h \in X_h$ to a local discrete solution $u_h \in X_h$ to $N_h(u_h)=0$ for \eqref{eqn:discrete} near an exact regular root $u \in X$ to \eqref{eqn:cont}, there is an algebraic error $\|u_h-v_h\|_{X_h}$ and an inconsistency error $\|v_h-Jv_h\|_{\widehat X}$ plus some intermediate residual $\rho(M)$.
Theorem~\ref{thm:apost} provides the equivalence of the error $\|u-v_h\|_{\widehat X}$ to
$$ \rho(M) + (1+M) \|v_h - Jv_h\|_{\widehat{X}} +M \|u_h-v_h\|_{X_h}$$
for any parameter $M \ge 0$. The underlying assumptions are phrased in a fairly general non-symmetric setting for rather general trilinear forms $\Gamma$ resp. $\Gamma_{\pw}$, and involve 
smallness assumptions on $\|u-u_h\|_{\widehat{X}} \le \varepsilon $ and $\|u_h-v_h\|_{X_h} \le \varrho$ that are guaranteed in the two applications to the stream-function vorticity formulation of the incompressible 2D Navier-Stokes and the  \vket~\cite{CCNNGRDS22}.
The provided a posteriori error analysis is generic and allows generalisations to other semilinear equations in future work.

\medskip \noindent Section~\ref{subsection:new} concerns the particular situation with $V=X=Y=H^2_0(\Omega)$ and {$V_h= X_h=Y_h \subseteq P_2(\T)$} for some triangulation $\T$ and a discrete norm $\|\bullet\|_h$ in $H^2(\T)$.
The arbitrary parameter $M$ in Theorem~\ref{thm:apost} becomes an upper bound of an interpolation operator $I_h: V \rightarrow V_h$ for the equivalence of $\rho(M)$  to an explicit residual-based a posteriori error estimator $\eta(\T) +\mu(\T)$ up to oscillations ${\rm osc}_k(F,\T)$.
The abstract parts of this paper in Section~\ref{sec:Abstract a posteriori error control} concludes with some remarks on the algebraic error $\|u_h-v_h\|_{V_h}$ in the context of the Newton-Kantorovich theorem and Section~\ref{subsection:new} illustrates the generality of the abstract results in Subsection~\ref{sec:algebraic}. In fact, given any approximation $v$ to the regular root $u$ of \eqref{eqn:cont}, that is a piecewise smooth function with respect to a triangulation $\T$ with $V_h \subset P_2(\T)$, the general Morley interpolation $I_\M$ and an adoption to $V_h$ (in case of $C^0$IP) leads to a postprocessed $v_h \in V_h$ that is close to $u$. If $\|u-v_h\|_{\widehat V}$ is sufficiently small, then \eqref{eqn:discrete} with $Q=S={J I_\M}$ provides a reference scheme such that Theorem~\ref{thm:apost} and \ref{thm:releff} provide a reliable and efficient estimate for $u-v_h$ (cf. Subsection~{2.5.2} for details).
Thus the a posteriori error analysis in this paper covers many more examples beyond the mandatory inexact solve in \eqref{eqn:discrete} or the computation of the discrete solution by a related scheme (e.g., without smoother for $R=S=Q={\rm id}$ in \eqref{eqn:discrete}).

\smallskip \noindent The application to the stream-function vorticity formulation of the incompressible 2D Navier-Stokes in Section~\ref{sec:nse}  enables the {\it first} explicit reliable and efficient residual-based a posteriori error estimates in this context and overcomes the gaps in \cite{CCGMNN_semilinear, kim_morley_2021}. The application to the \vket~in Section~\ref{sec:vke} also considers single forces in the source terms. %
Section~\ref{sec:Numerical experiments} presents the {\it first} numerical comparisons of the quadratic schemes and confirms the a priori equivalence results in \cite{CCNNGRDS22}. The associated adaptive mesh-refining recovers the optimal convergence rates.

\smallskip \noindent
Standard notation on  Lebesgue and Sobolev spaces, 
their norms, and $L^2$ scalar products applies throughout this paper;  $\|\bullet\|$ %
abbreviates the operator norm of a linear operator.  
The Hilbert space $V := H_0^2(\Omega)$ is endowed with the energy scalar product $a(\bullet,
\bullet)$ that induces
the $H^2$ seminorm $\displaystyle\trinl\bullet\trinr := |\bullet|_{H^2(\Omega)}$; the induced dual linear operator norm in $H^{-2}(\Omega)$ is denoted by $\trinl \bullet \trinr_{*}$ in the later sections. In the sequel, the notation $A \lesssim B$ abbreviates $A \leq CB$ for some positive generic constant $C$, 
which exclusively depends on the shape-regularity of the underlying triangulation $\T$ (i.e., on $0 < \omega_0 \le {\arcangle} {T}$),  
 $A\approx B$ abbreviates $A\lesssim B \lesssim A$. 

\section{Abstract a posteriori error analysis}\label{sec:Abstract a posteriori error control}
The a posteriori error analysis concerns some approximation $v_h \in X_h$ in some
discrete nonconforming space $X_h\not\subset X$ to a regular root $u$ of the continuous problem \eqref{eqn:cont}. The known approximation could result from an inexact solve of the discrete problem 
\eqref{eqn:discrete} with local solution $u_h \in X_h$; hence $u_h$ and $v_h$ are different in general and the main interest is on the distance of $u$ and $v_h$. The abstract results of this section also apply to semilinear second-order problems in future work. 
		
		\subsection{Discretisation}%
		\label{sub:Discretisation}

\noindent Let $\widehat{X}$ (resp. $\widehat{Y}$) be a real Banach space with norm $\|\bullet\|_{\widehat{X}}$ (resp.\ $\|\bullet\|_{\widehat{Y}}$) and suppose $X$ and $X_h$ (resp.\ $Y$ and $Y_h$) are two complete linear subspaces of $\widehat{X}$ (resp.\ $\widehat{Y}$) with inherited norms $\|\bullet\|_{X}:=\big(\|\bullet\|_{\widehat{X}}\big)|_{X}$ and $\|\bullet\|_{X_h}:=\big(\|\bullet\|_{\widehat{X}}\big)|_{X_h}$ (resp.\ $\|\bullet\|_{Y}:=\big(\|\bullet\|_{\widehat{Y}}\big)|_{Y}$ and $\|\bullet\|_{Y_h}:=\big(\|\bullet\|_{\widehat{Y}}\big)|_{Y_h}$){; $X + X_h \subseteq \widehat{X}$ and $Y + Y_h \subseteq \widehat{Y}$.} Let the bounded linear operator
${A}\in L({X}; {Y}^*)$ be associated to the 
bilinear form $a$ and suppose $A$ is invertible and, in particular, satisfies
\begin{gather} 
0 < \alpha :=\inf_{\substack{{x} \in {X} \\ \|{x}\|_{{X}}=1}} \sup_{\substack{{y}\in {Y}\\ \|{y}\|_{{Y}}=1}} a({x},{y}).
\label{cts_infsup}
\end{gather}
Let $ \Gamma_{\pw}:\widehat{X}\times \widehat{X}\times \widehat{Y}\to \R$
denote a bounded trilinear form that extends $\Gamma=\Gamma_{\pw}|_{X\times X\times Y}$ 
such that
$$\displaystyle \| \Gamma_{\pw}\|:=\| \Gamma_{\pw}\|_{\widehat{X}\times \widehat{X}\times \widehat{Y}}:=
\sup_{\substack{\widehat{x}\in \widehat{X}\\ \|\widehat{x}\|_{\widehat{X}}=1}}
\sup_{\substack{\widehat{\xi}\in \widehat{X}\\ \|\widehat{\xi}\|_{\widehat{X}}=1}}
\sup_{\substack{\widehat{y}\in \widehat{Y}\\ \|\widehat{y}\|_{\widehat{Y}}=1}}
\Gamma_{\pw}(\widehat{x},\widehat{\xi},\widehat{y})~<\infty\quad\text{and set}\quad\| \Gamma\|:=\| \Gamma\|_{{X}\times {X}\times {Y}}.$$
Define the quadratic function $N:X \rightarrow Y^*$ by \eqref{eqn:cont}.
A vector $u \in X$ is called a {\it regular root} of \eqref{eqn:cont},  if  $u$ solves $N(u)=0$
 and the Frech\'et derivative $DN(u)$ is a bijection and, in particular, fulfils 
\begin{gather} 
0<\beta:=\inf_{\substack{x\in X\\ \|x\|_{X}=1}} \sup_{\substack{y\in Y\\ \|y\|_{Y}=1}}\Big(a(x,y)+ {\Gamma}(u,x,y)+{\Gamma}(x,u,y)\Big).
\label{cts_infsup_b}
\end{gather}
Suppose that the bounded bilinear form $a_h:X_h \times Y_h \rightarrow {\mathbb R}$ suffices the discrete inf-sup condition
\begin{align} \label{eqn:stability}
&0 <  \alpha_h  :=\inf_{\substack{{x}_h\in {X}_h \\ \|{x_h}\|_{{X_h}}=1}} \sup_{\substack{{y}\in {Y_h}\\ \|{y_h}\|_{{Y}_h}=1}}\dislin({x_h},{y_h})
\end{align}
for some constant $\alpha_h$. 
The quasi-optimal smoothers are linear and bounded operators  $P \in L(X_h;X)$, $Q \in L(Y_h;Y)$,  $\rop \in L(X_h;\widehat{X})$, $S \in L(Y_h;
\widehat{Y})$ {with respective operator norms $\|P\|, \|Q\|, \|R\|$, and $\|S\|$}  such that, for all $x_h \in X_h, \: x \in X,\: y_h \in Y_h, $ and $y \in Y$,
\begin{align}
	\| (1 - P)x_h \|_{\widehat{X}} & \le \Lambda_{\rm P} \| x - x_h \|_{\widehat{X}}, \label{quasioptimalsmootherP}\\
	\| (1 - Q)y_h \|_{\widehat{Y}} & \le \Lambda_{\rm Q} \| y - y_h \|_{\widehat{Y}}, \label{quasioptimalsmootherQ}\\
	\| (1 - \rop)x_h \|_{\widehat{X}} & \le \lamr\| x - x_h \|_{\widehat{X}}, \label{quasioptimalsmootherR}\\
	\| (1 - \sop)y_h \|_{\widehat{Y}} & \le \lams\| y - y_h \|_{\widehat{Y}} \label{quasioptimalsmootherS}
 \end{align}
hold with constants $\Lambda_{\rm P},\Lambda_{\rm Q},  \lamr, \lams \ge 0$. Suppose there exists $\Lambda_{\rm C}>0$ such that
\begin{align}
 a(Px_h,Qy_h)- \dislin(x_h,y_h) & \le \Lambda_{\rm C} \| x_h - Px_h\|_{\widehat{X}} \| y_h \|_{Y_h} \label{eqn:H}
\end{align}
holds for all $(x_h,y_h) \in X_h \times Y_h$. While \eqref{eqn:stability} is stability, \eqref{eqn:H} is consistency introduced in \cite{ccnnlower2022} for linear problems. Let the quadratic function $N_h:X_h \rightarrow Y_h^*$ be defined by \eqref{eqn:discrete}. The local conditions on the roots $u$, $u_h$, and their approximation $v_h$ are summarised as follows. 
\begin{enumerate}
\item[{\bf (L)}]
 Let $u\in X$ denote a regular root of \eqref{eqn:cont} for a given source term $F \in  Y^*$ and
    let there exist $\varepsilon,\varrho >0$ and  {$0<\kappa <1$} such that
    \begin{itemize}
    \item[{\bf (L1)}]  $N_h(u_h)=0$  holds for exactly one solution $u_h\in X_h $  with $\| u-u_h\|_{\widehat{X}} \le \varepsilon$, 
    \item[{\bf (L2)}] $v_h \in X_h$ satisfies $\|u_h-v_h\|_{{X_h}}  \le \varrho$,
    \item[{\bf (L3)}] $\varepsilon+\varrho\leq\kappa\beta/\left((1+\LP)\|\Gamma\|\right)$.
    \end{itemize}
    \end{enumerate}
The point is that the recent paper \cite{CCNNGRDS22} provides affirmative examples for all those conditions \eqref{cts_infsup}-\eqref{eqn:H}, {\bf (L1)}-{\bf (L3)}, and \eqref{eqn:QO} with constants that are independent of some discretisation parameter {$h$}, provided the discretization is sufficiently fine.
\subsection{Abstract a posteriori analysis} 
\label{subsec:outline}

This section presents an  abstract reliability  and  efficiency result. %
The abstract a posteriori error control has three contributions. The first one is an  intermediate residual 
\begin{align} \label{a}
\displaystyle \rho(M)&:= \sup_{\stackrel{y \in Y}{\|y \|_Y \le 1} } \inf_{\stackrel{y_h \in Y_h}{\|y_h \|_{Y_h} \le M} }  %
\Big( F(y-Qy_h)-a(Pv_h,y-Qy_h) -
\Gamma_{\pw}(Rv_h, Rv_h, y-Sy_h)  \Big)
\end{align}
for some parameter $M \ge 0$. The role of $M$ will be clarified in Section~\ref{subsection:new} below. At this point it suffices to observe that $y_h$ may be some discrete object such that $Qy_h$ and $Sy_h$ approximate $y$ and we expect
$\|y_h\|_{Y_h} \le M \|y\|_{Y} \le M$ is bounded.  Notice that $\rho$ is monotone decreasing and $\rho(0) \ge \rho(M)$. %
The second contribution is a consistency term $\|v_h-Pv_h\|_{\widehat{X} }$ (computable from $v_h \in X_h$ and
the quasi-optimal smoother $P$) with \eqref{quasioptimalsmootherP} and throughout serves as an efficient a posteriori term
\cite{CCBGNN22}. The third term $\|u_h-v_h\|_{X_h} \le \varrho$ is the algebraic error (e.g., from an inexact solve) and  is briefly discussed in Subsection~\ref{sec:algebraic}.
\begin{theorem}[abstract reliability and efficiency]\label{thm:apost}
	(a) Suppose {\bf (L1)}--{\bf (L3)}, \eqref{quasioptimalsmootherP}-\eqref{quasioptimalsmootherR}, \eqref{eqn:H}, and $M \ge 0$. Then 
\begin{align}
  & \displaystyle \| u- v_h\|_{\widehat{X}}  \le
C_{\rm rel} \Big( \rho(M)  + (1+M) \|(1-P)v_h\|_{\widehat{X}} +M\|u_h-v_h\|_{X_h} \Big). %
 \label{eqn:apost} %
\end{align}
The  constant  $C_{\rm rel} $  %
exclusively depends on  $\beta, (1-\kappa)^{-1}, \Lambda_{\rm C}, \Lambda_{\rm P},$ $ \Lambda_{\rm R}, \|a_h\|, $ $ \|\Gamma\|,\|\Gamma_\pw\|,\|S\|,$ and $\|u\|_{X}$.

\medskip \noindent 
(b) Suppose {\bf (L1)}--{\bf (L2)} and \eqref{eqn:QO}. Then $\rho(M) \le \rho(0)$ and
\begin{align} \label{eq:eff}
& \rho(0)  +\|(1-P)v_h\|_{\widehat{X}} +\|u_h-v_h\|_{X_h} 
 \le C_{\rm eff}\| u- v_h\|_{\widehat{X}}. 
 \end{align}
 The  constant  $C_{\rm eff} $  %
exclusively depends on $ \beta, (1-\kappa)^{-1}, \Lambda_{\rm C},\Lambda_{\rm P}, \Lambda_{\rm R}, \|a_h\|,  \|\Gamma\|,\|\Gamma_\pw\|,\|S\|,\|u\|_{X}$, and $C_{\rm qo}$.
\end{theorem}

\subsection{Proof of Theorem~\ref{thm:apost}.a}
	\label{sec:apost}
The proof is split into several subsections below with $u$, $u_h$, $v_h$ as in the statement of Theorem~\ref{thm:apost}.
\subsubsection{Reduction to $\|u-Pv_h\|_X$}\label{step1}
Consequences of  \eqref{quasioptimalsmootherP}, \eqref{quasioptimalsmootherR},  and a  triangle inequality read, for 
 any
$x\in X$ and $x_h \in X_h$, as
\begin{align} 
\|v_h-Pv_h\|_{\widehat{X}} &\le  \Lambda_{\rm P} \|u-v_h\|_{\widehat{X}}, \label{eqn:vh-Pvh} \\
 \|x- Px_h\|_{{X}} &\le (1+\Lambda_{\rm P}) \|x-x_h\|_{\widehat{X}},  \label{eqn:PR1} \\
\|x- Rx_h\|_{\widehat{X}} &\le (1+\Lambda_{\rm R}) \|x-x_h\|_{\widehat{X}},  \label{eqn:PR2} \\
 \| (P-R)x_h\|_{\widehat{X}} & \le (1+\Lambda_{\rm R}) \|  x_h- Px_h\|_{\widehat{X}}. \label{eqn:PR3}
\end{align}
The efficiency \eqref{eqn:vh-Pvh} of the  a posteriori estimator $\|v_h-Pv_h\|_{\widehat{X}}$ and a triangle inequality 
	\begin{equation}\label{eq.tri.rel}
	\|u-v_h\|_{\widehat{X}} \le \|u-Pv_h\|_{{X}}+\|v_h-Pv_h\|_{\widehat{X}} 
	\end{equation}
motivate the focus on $\|u-Pv_h\|_{{X}}$ in the error analysis below.
\subsubsection{Reduction to $ \| N(Pv_h)\|_{Y^*}$}
	The inf-sup condition \eqref{cts_infsup_b} with $\beta>0$  for the regular root $u$ %
 leads, for any $\tau>0$, to
some $y\in Y$ with $\|y\|_{Y}\le 1+\tau$  and 
	\begin{equation}\label{Der_infsup}
	\beta \|u-Pv_h\|_{X} =  DN(u;Pv_h-u,y).
	\end{equation}
 (For reflexive Banach spaces, $\tau=0$ is possible, but for the time being we require $\tau>0$). 
	Since $N$ is quadratic, the finite Taylor series is exact, namely
	\begin{align*}%
	N(Pv_h;y)=N(u;y)+DN(u;Pv_h-u,y)+\half D^2N(u;u-Pv_h,u-Pv_h,y).
	\end{align*}
	Since $N(u)=0$ and $D^2N(u;u-Pv_h,u-Pv_h,y)=2 \, \Gamma(u-Pv_h,u-Pv_h,y)$, this reads 
	\begin{align}\label{quadratic_identity}
	N(Pv_h;y)+DN(u;u-Pv_h,y)& =  \Gamma(u-Pv_h,u-Pv_h,y).
	\end{align}
	The combination of \eqref{Der_infsup}-\eqref{quadratic_identity} and the bound $\|\Gamma\|$ of the trilinear form result in 
	\begin{align}\label{eqn:newa}
{\beta}  \|u-Pv_h\|_{X}& %
	\leq \big(\| N(Pv_h)\|_{Y^*} + \|\Gamma\| \|u-Pv_h\|^2_{X}\big) (1+\tau)
	\end{align}
	with $\|y\|_{Y} \le 1+\tau$ in the last step. 
Recall that \eqref{eqn:newa} holds for any $\tau>0$ and so $\tau \searrow0 $ provides
\begin{align} \label{eqn:new}
{\beta}  \|u-Pv_h\|_{X}& \le \| N(Pv_h)\|_{Y^*} + \|\Gamma\| \|u-Pv_h\|^2_{X}.
\end{align}
Since $\|u-Pv_h\|_{X} \le (1+\Lambda_{\rm P}) \|u-v_h\|_{\widehat{X}}$ by \eqref{eqn:PR1} and
$\|u-v_h\|_{\widehat{X}} \le \|u-u_h\|_{\widehat{X}} + 
\|u_h-v_h\|_{X_h} \le \varepsilon + \varrho$ by {\bf (L1)}-{\bf (L2)}, we infer
$$ \|\Gamma \| \|u-Pv_h\|_{X} \le 
 \|\Gamma \| (1+\Lambda_{\rm P}) (\varepsilon + \varrho) \le \kappa \beta$$ with  {\bf (L3)} in the last step. This and \eqref{eqn:new} imply
 \begin{align} \label{eqn:bounda}
 \beta (1-\kappa) \|u-Pv_h\|_{X}
 & \le \|N(Pv_h)\|_{Y^*}.
 \end{align}
 The combination of \eqref{eq.tri.rel} and \eqref{eqn:bounda} reveals
 \begin{align} \label{eqn:u-vh}
    \|u-v_h\|_{\widehat{X}} & \le 
     \|v_h -Pv_h\|_{\widehat{X}}+
    \beta^{-1} (1-\kappa)^{-1} \|N(Pv_h)\|_{Y^*}
 \end{align}
 and we are left with the a posteriori analysis of $\|N(Pv_h)\|_{Y^*}$.
\subsubsection{Appearance of $\rho(M)$}\label{sub:varrho} To control $\| N(Pv_h)\|_{Y^*}$, consider any $y \in Y$ with $\|y\|_{Y} =1$
and any $y_h \in Y_h$ with $\|y_h\|_{Y_h} \le M$. 
Elementary algebra with the definition of $N(Pv_h; y)$ leads to
\begin{align} \label{eqn:newb}
	N(Pv_h; y)& = a(Pv_h,y-Qy_h) + \Gamma_{\pw}(Rv_h,Rv_h, y-Sy_h) -{F}(y-Qy_h) \nonumber \\
	& \quad + \Gamma_{\pw}(Pv_h,Pv_h,y) -\Gamma_{\pw}(Rv_h,Rv_h, y)  
	\nonumber \\
	& \quad  +a(Pv_h,Qy_h)-F(Qy_h) + \Gamma_{\pw}(Rv_h,Rv_h, Sy_h) =:S_1+S_2+S_3.
\end{align}
The  first  term  $S_1$ %
 gives rise to the intermediate residual 
$$
S_1:=a(Pv_h,y-Qy_h) + \Gamma_{\pw}(Rv_h,Rv_h, y-Sy_h) -{F}(y-Qy_h) \le \rho(M) $$
provided $y_h \in Y_h$ is selected to obtain an infimum in \eqref{a}. \big(The analysis of the remaining terms $ S_2+S_3 $ exclusively utilizes  $y_h \in Y_h$ and $\|y_h\|_{Y_h} \le M$  below.\big) Thus
\begin{align}\label{eqn:intermediate}
N(Pv_h; y)
 \le  \rho(M) + S_2+S_3. 
 \end{align}
\subsubsection{Difference of the trilinear form $S_2$}\label{S2}
Elementary algebra and the boundedness  of the piecewise trilinear form result in 
\begin{align} \label{t2}
S_2 & :=\Gamma_{\pw}(Pv_h,Pv_h,y) -\Gamma_{\pw}(Rv_h,Rv_h, y)   \nonumber \\
&= \Gamma_{\pw}((P-R)v_h,Pv_h,y) +\Gamma_{\pw}(Rv_h,(P-R)v_h, y) \nonumber \\
& \le \|\Gamma_\pw\|  (1+\Lambda_{\rm R}) (\|P\|+ \|R\|) \|v_h\|_{X_h} \|v_h-Pv_h \|_{\widehat{X}}
\end{align}
with \eqref{eqn:PR3} and $\| y\|_Y {=1}$  in the last step.

\medskip \noindent 
All the operator norms $\|P\|, \|Q\|, \|R\|, \|S\|$ of the quasi-optimal smoothers are controlled in terms of  $\Lambda_{\rm P}$, $\Lambda_{\rm Q}$,  $\Lambda_{\rm R}$, and  $\Lambda_{\rm S}$. For instance, \eqref{quasioptimalsmootherR} shows  $\|Rv_h\|_{\widehat{X}} \le (1+\Lambda_{\rm R})  \|v_h\|_{X_h}$.
  \subsubsection{Remaining bound $S_3$}
The last term $ S_3$ on the right-hand side of \eqref{eqn:newb} reads 
$$S_3:= a(Pv_h,Qy_h)-F(Qy_h) + \Gamma_{\pw}(Rv_h,Rv_h, Sy_h).$$
A comparison with
$N_h(u_h;y_h)=0$ from \eqref{eqn:discrete} and elementary algebra result in 
\begin{align}\label{eqn:s3}
 S_3&= a(Pv_h,Qy_h)- a_h(v_h,y_h ) + a_h(v_h-u_h,y_h) \nonumber \\
 & \qquad + \Gamma_{\pw}(Rv_h,Rv_h, Sy_h) - \Gamma_{\pw}(Ru_h,Ru_h, Sy_h).
 \end{align}
The consistency \eqref{eqn:H} controls the first two terms in the right-hand side of \eqref{eqn:s3},
\begin{align}\label{eqn:H_application}
a(Pv_h,Qy_h)- a_h(v_h,y_h ) & \le \Lambda_{\rm C} \|v_h-Pv_h\|_{\widehat{X}} \|y_h\|_{Y_h}.
\end{align}
The boundedness of $a_h$ establishes 
$a_h(v_h-u_h,y_h) \le \|a_h\|\|v_h-u_h\|_{X_h} \|y_h\|_{Y_h}.$
Elementary algebra  for the last two terms in \eqref{eqn:s3} provides
\begin{align} \label{eqn:t2a}
& \Gamma_{\pw}(Rv_h,Rv_h, Sy_h) - \Gamma_{\pw}(Ru_h,Ru_h, Sy_h)  \nonumber \\
& =\Gamma_{\pw}(R(v_h-u_h),Rv_h, Sy_h)  + \Gamma_{\pw}(Ru_h,R(v_h-u_h), Sy_h)  \nonumber \\
 &  \le 
  \|\Gamma_\pw\|   \|R\|^2 \|S\| \Big( \|v_h\|_{X_h} +\|u_h\|_{X_h} \Big)  
   \|u_h-v_h\|_{X_h} \|y_h\|_{Y_h}
  \end{align}
  with boundedness of the piecewise trilinear form in the last step. 
A combination of the aforementioned  estimates with \eqref{eqn:s3} and $\|y_h\|_{Y_h} \le M$  shows
\begin{align}\label{eqn:s3a}
S_3 & \le M \Lambda_{\rm C} \|v_h-Pv_h\|_{\widehat{X}} + M\Big( \|a_h\|  + \|\Gamma_\pw\|   \|R\|^2 \|S\|  \big( \|v_h\|_{X_h} +\|u_h\|_{X_h} \big)  \Big)  \|u_h-v_h\|_{X_h}.
\end{align}
\subsubsection{Final a posteriori error estimate}\label{ssub:abstract_reliability_final}
Since $y \in Y$ with $\|y\|_Y =1$ is arbitrary, the combination of \eqref{t2} and \eqref{eqn:s3a} in \eqref{eqn:intermediate} leads to 
\begin{align}
\|N(Pv_h)\|_{Y^*} & \le \rho(M) + \Big(M \Lambda_{\rm C}  + \|\Gamma_\pw\|  \left(1+\Lambda_{\rm R} \right) \left(\|P\|+ \|R\|\right) \|v_h\|_{X_h}\Big) \|v_h-Pv_h\|_{\widehat{X}} \nonumber \\
& \qquad + M \Big( 
\|a_h\| + \|\Gamma_\pw\|   \|R\|^2 \|S\| 
\big( \|v_h\|_{X_h} +\|u_h\|_{X_h} \big)
 \Big) \|u_h-v_h\|_{X_h}.  \label{eqn:npvh}
\end{align}
Triangle inequalities and {\bf (L1)}-{\bf (L3)} reveal 
\begin{align} \label{eqn:vhbound}
\|v_h\|_{X_h} &\le \|u\|_X + \|u-u_h\|_{\widehat{X}} +\|u_h-v_h\|_{X_h} \le\|u\|_{X} + \varepsilon + \varrho  \le \|u\|_{X} +  \beta/ ((1+\Lambda_{\rm P}) \|\Gamma\|).
\end{align}
The same arguments apply to show $\|u_h\|_{X_h} \le \|u\|_X + \|u-u_h\|_{\widehat{X}} \le \|u\|_{X} +  \beta/ ((1+\Lambda_{\rm P}) \|\Gamma\|)$.
 A substitution of \eqref{eqn:vhbound} and the analog estimate for $\|u_h\|_{X_h}$ in 
\eqref{eqn:npvh} reveal
\begin{align}\label{eqn:ta}
\|N(Pv_h)\|_{Y^*}  \le \rho(M) + (M\Lambda_{\rm C} + C_1)  \|v_h-Pv_h\|_{\widehat{X}}  +C_2 M \|u_h-v_h\|_{X_h}
\end{align}
with  universal constants $\constant{}\label{ct4}:= 
\|\Gamma_\pw\|  \left(1+\Lambda_{\rm R} \right) \left(\|P\|+ \|R\|\right) \big(\|u\|_X+ \beta/ \big((1+\Lambda_{\rm P}) \|\Gamma\| \big) \big)$ and $\constant{} \label{ct5}:=
 \|a_h\| + 2\|\Gamma_\pw\|   \|R\|^2 \|S\|   \big(\|u\|_X+ \beta/ \big((1+\Lambda_{\rm P}) \|\Gamma\| \big) \big)
 $.
A combination \eqref{eqn:ta} with \eqref{eqn:u-vh} provides
\begin{align} \label{eqn:est1}
\beta(1-\kappa) \|u-v_h\|_{\widehat{X}
 }& \le   \rho(M) + C_3 (1+M)\|v_h-Pv_h\|_{\widehat{X}}  +C_2 M\|u_h-v_h\|_{X_h}
\end{align}
with $\constant{}\label{ct3} := \max \big{\{} \Lambda_{\rm C}, C_1+\beta(1-\kappa) \big{\}} $. This concludes the proof of reliability with a reliability constant 
$C_{\rm rel}:= \beta^{-1}(1-\kappa)^{-1}\max \big{\{} 1, C_2 M, C_3 (1+M)\big{\}} $. \qed
\subsection{Proof of Theorem~\ref{thm:apost}.b}\label{sec:eff}
The efficiency of $\|v_h-Pv_h\|_{\widehat{X}}$ 
follows from \eqref{eqn:vh-Pvh} and hence the focus is on the other two terms $\rho(M)\leq\rho(0)$ and $\|u_h-v_h\|_{X_h}$.
\subsubsection{Linear intermediate problem}
The link between $\rho(M)$ from \eqref{a} to the known a posteriori results for linear problems reviewed in \cite{CCBGNN22} is the linear intermediate problem
\begin{align}\label{c}
a(\widetilde{u},y)&=F(y) - \Gamma_{\pw}(Rv_h,Rv_h,y)  \text{  for all } y \in Y.
\end{align}
Since the associated operator $A:X\to Y^*$ is invertible, the problem \eqref{c} admits a unique
solution $\widetilde u\in X$.
It follows that 
\begin{align}\label{d}
 \rho(M) & \le \rho(0)=   \sup_{\stackrel{y \in Y}{\|y \|_Y \le 1} } a(\widetilde{u}-Pv_h,y) \le \|a\| \: \|\widetilde{u} -Pv_h\|_{X}.
 \end{align} 
\subsubsection{Efficiency of $ \|u-\widetilde{u} \|_{X}$}
Recall that $u$ is a fixed regular root of $N$, while $\widetilde{u}$ solves 
\eqref{c}.  %
The inf-sup condition in \eqref{cts_infsup} leads,  for any $\tau>0$, to some  $ y \in Y$  with $\|y\|_Y \le 1+\tau$ and 
\begin{align*} 
\alpha \|u-\widetilde{u}\|_X &= a(\widetilde{u}-u,y) = \Gamma(u,u,y) -\Gamma_{\pw}(Rv_h,Rv_h,y )%
\end{align*}
with \eqref{eqn:cont}  and \eqref{c} in the last step. This, the boundedness of the trilinear form, $\|u-Rv_h\|_{\widehat{X}} \le (1+\Lambda_{\rm R}) \|u-v_h\|_{\widehat{X}}$ from  \eqref{eqn:PR2}, and $\|y\|_Y \le 1+\tau$  provide
\begin{align*} 
\alpha \|u-\widetilde{u}\|_X &=\Gamma_{\pw}(u-Rv_h,u,y ) + \Gamma_{\pw}(Rv_h,u-Rv_h,y ) \nonumber \\
& 
 \le  \|\Gamma_\pw\| (1+\tau) (1+\Lambda_{\rm R}) \|u -v_h\|_{\widehat{X}} (\|u\|_X + \|Rv_h\|_{\widehat{X}}).
\end{align*}
The aforementioned estimate holds for any $\tau>0$, hence $\tau \searrow0 $ and  %
\eqref{eqn:vhbound} establish
\begin{equation} \label{eqn:u-ut}
 \|u- \widetilde{u}\|_{X} \le %
 \alpha^{-1} \|\Gamma_\pw\|  (1+\Lambda_{\rm R}) \Big(\|u\|_X + \|R\| \big(\|u\|_X +  \beta/ ((1+\Lambda_{\rm P}) \|\Gamma\|) \big)   \Big) \|u-v_h\|_{\widehat{X}}.
 \end{equation}
\subsubsection{Efficiency of $\rho(M)$}The intermediate problem \eqref{c} leads to 
\eqref{d}, namely 
\begin{align}\label{a1}
\rho(M) \le\rho(0) %
& \le \| a\| \| \widetilde{u} -Pv_h\|_{\widehat{X}} \le \|a\| \left(\|u - \widetilde{u}\|_X+ \|u-Pv_h\|_{X}\right).
\end{align}
Recall $\| u -Pv_h\|_{\widehat{X}} \le (1+\Lambda_{\rm P}) \| u -v_h\|_{\widehat{X}}$ from \eqref{eqn:PR1} and combine it with  \eqref{eqn:u-ut}-\eqref{a1} to deduce
$$\rho(0)  \le \|a\|\bigg(1+ \Lambda_{\rm P} +\alpha^{-1} \|\Gamma_\pw\|  (1+\Lambda_{\rm R}) \Big(\|u\|_X + \|R\| \big(\|u\|_X +  \beta/ ((1+\Lambda_{\rm P}) \|\Gamma\|) \big)   \Big)\bigg) \|u-v_h\|_{\widehat{X}}. \qed$$
\subsubsection
{Efficiency of $\|u_h-v_h\|_{X_h}$ under \eqref{eqn:QO}}  The quasi-best approximation \eqref{eqn:QO} implies $\|u-u_h\|_{\widehat{X}} \le C_{\rm qo} \|u-v_h\|_{\widehat{X}}$. This and a triangle inequality provide
$ \|u_h-v_h\|_{X_h} \le \|u-u_h\|_{\widehat{X}} + \|u-v_h\|_{\widehat{X}} \le (1+C_{\rm qo})\|u-v_h\|_{\widehat{X}}.$ \qed

\subsection{Comments}

\subsubsection{Inexact solve}
The ad hoc application of Theorem~\ref{thm:apost} is on \eqref{eqn:discrete} with $S=Q$ and a flexible choice of $R$ with \eqref{quasioptimalsmootherR}. The local convergence of the Newton scheme is guaranteed in \cite{CCNNGRDS22} and a termination leads to $v_h \in X_h$ with an algebraic error discussed in Subsection~\ref{sec:algebraic} below. 
A few iterations more 
 provide the discrete solution up to machine precision and $\|u_h-v_h\|_h$ is negligible and this point of view is adapted in Section~\ref{sec:Numerical experiments}. %
\subsubsection{Approximation $v_h$ from other discretisations}
\label{sub:Approximation from other discretisations}
The discrete scheme in \eqref{eqn:discrete} with smoother $S=Q$ for the definition of $N_h$ in \eqref{eqn:discrete} and its root $u_h \in X_h$ with $\|u-u_h\|_{\widehat{X}} \le \varepsilon$ from {\bf(L1)} can serve as a reference scheme. Given an accurate approximation $v_h  \in X_h$ from another  numerical scheme that is sufficiently good in the sense that
\begin{equation}\label{eqn:varepsilon}
\|u-v_h\|_{\widehat{X}} \le \min \left\{\varepsilon,  {\varrho}/{(1+ C_{\rm qo})} \right\}.
\end{equation}
Since \eqref{eqn:QO} provides $\|u-u_h\|_{\widehat{X}} \le C_{\rm qo} \|u-v_h\|_{\widehat{X}}$, a triangle inequality and \eqref{eqn:varepsilon} reveal
\[  \| u_h-v_h\|_{X_h} \le (1+C_{\rm qo}) \| u-v_h\|_{\widehat{X}} \le \varrho.\]
Hence Theorem~\ref{thm:apost} applies to $v_h \in X_h$ and the explicit residual-based a posteriori estimators $\eta + \mu$ of Subsection~\ref{subsection:new}  lead to reliable and efficient error control of $\| u-v_h\|_{\widehat{X}} $. Although $v_h$ may originate from a very different setting, its a posteriori error control, namely the evaluation of $\| u-v_h\|_{\widehat{X}}$ in Theorem~\ref{thm:apost}, concerns the reference scheme $N_h$ from \eqref{eqn:discrete} with  $S=Q$.

\subsubsection{Control of algebraic errors}\label{sec:algebraic}
{The numerical analysis of the discrete problem as a high-dimensional algebraic system of equations %
is a routine task, e.g., with the known Newton scheme and the Newton-Kantorovich theorem. For instance, suppose that 
 $DN(u_h)$ satisfies the discrete inf-sup condition
\begin{align} \label{eqn:discreteinfsup}
0< \beta_h:=\inf_{\substack{x_h\in X_h\\ \|x_h\|_{X_h}=1}} \sup_{\substack{y_h\in Y_h\\ \|y_h\|_{Y}=1}}\Big(a_h(x_h,y_h)+ {\Gamma}_{\pw}(u_h,x_h,y_h)+{\Gamma}_{\pw}(x_h,u_h,y_h)\Big)
\end{align}
proved in  \cite{CCNNGRDS22} and recall that $u_h\in X_h$ is a discrete root of \eqref{eqn:discrete}. 
\begin{lemma}[control of $\|u_h-v_h\|_{X_h}$]\label{lem:uh-vh} Any $v_h\in X_h$ and $0<\kappa<1$ with $\|u_h-v_h\|_{X_h}\leq \kappa\beta_h/\big(\|\Gamma_\pw\|\|R\|^2\|S\|\big)$ satisfy
	\begin{align*}
		(1-\kappa)\beta_h\|u_h-v_h\|_{X_h} \leq \|N_h(v_h)\|_{Y_h^*}\leq \big(\kappa\beta_h +
		\|DN_h(u_h)\|_{X_h^*\times Y_h^*}\big)\|u_h-v_h\|_{X_h}.
	\end{align*}
\end{lemma}
\noindent The lemma is proved in Supplement A and an associated  termination criterion is outlined in Supplement~C.
\begin{example}[Computation of  $\|N_h(v_h)\|_{Y_h^*}$] \label{ex:computation_residual}If 
${\rm dim}(X_h)= {\rm dim}(Y_h) < \infty$ and 
\eqref{eqn:stability} holds, the linear operator $A_h : X_h \rightarrow Y_h^*$ associated with the bilinear form $a_h: X_h \times Y_h \rightarrow {\mathbb R}$ is invertible with  $\|A_h\|_{L(X_h; Y_h^*)}\|A_h^{-1}\|_{L(Y_h^*; X_h)} \le \|a_h\|/\alpha_h.$ Hence one linear solve of 
$a_h(\xi_h, y_h) =N_h(v_h; y_h)$ for a unique discrete solution $\xi_h \in X_h$ suffices for $\alpha_h \|\xi_h\|_{X_h} \le \|N_h(v_h)\|_{Y_h^*} \le \|a_h\| \|\xi_h\|_{X_h}$ and makes Lemma~\ref{lem:uh-vh} applicable.
\end{example}
\section{Explicit residual-based a posteriori estimator}\label{subsection:new}
This section discusses computable and explicit bounds for the intermediate residual $\rho(M)$ (and the consistency term $\|v_h-Pv_h\|_{X_h}$) in the reliablity control of Theorem \ref{thm:apost} in an application to fourth-order semilinear problems.
\subsection{Triangulation, interpolation, and smoother}\label{sub:triangulation, interpolation, and smoother}
Throughout this paper, $\T$ denotes a shape-regular triangulation of a polygonal and bounded (possibly
multiply-connected) Lipschitz domain $\Omega\subset \mathbb
R^2$ into
triangles. The set of all  vertices $\mathcal{V}$ (resp.\ edges $\E$)  in the triangulation $\cT$  decomposes into interior vertices  $\mathcal{V}(\Omega)$ (resp.\ interior edges $\cE(\Omega)$) and boundary vertices  $\mathcal{V} (\partial\Omega)$ (resp.\ boundary edges $\cE(\partial\Omega)$).
Let $h_E\coloneqq |E|\coloneqq\mathrm{diam}(E)=|A-B|$ denote the length of any edge $E=\mathrm{conv}(A,B)\in\E$ with vertices $\mathcal V(E)=\{A,B\}$.
Define the piecewise constant mesh size $h_{\cT}(x)=h_T={\rm diam}  (T)$ for all $x \in T\in \cT$ (resp. $h_{\cE}(x)=h_E={\rm diam}  (E)$ for all $x \in E \in \cE$), and set $h_{\rm max} :=\max_{T\in \cT}h_T$.
The notation $\bT(\delta)$ denotes a family of  those triangulations $\T$ with maximal-mesh size $h_{\rm max}\leq\delta$ smaller than or equal to $\delta>0$ and interior angles $\ge \omega_0 >0$ for some universal constant $\omega_0$. 

The space $P_{\hspace{-.13em}k}(T)$ of polynomials of total degree at most $k\in\mathbb N_0$ on $T\in\T$ defines
the space of piecewise polynomials
\begin{align*}
	P_{\hspace{-.13em}k}(\T)&\coloneqq\{p\in L^\infty(\Omega) : p|_{T}\in P_{\hspace{-.13em}k}(T)\text{ for all } T\in \T\}
\end{align*}
and let $\Pi_k$ denote the $L^2$ projection onto $P_k(\T)$; $\Pi_k$ acts componentwise on vectors or matrices. Here and throughout this paper, $H^m(\T)\coloneqq \prod_{T\in\T}H^m(T)$ is the space of piecewise Sobolev functions for $m=1,2$ with the abbreviation $H^m(K)\coloneqq
H^m(\mathrm{int}\;K)$ for a triangle or edge $K\in \T\cup \E$ with relative interior $\mathrm{int}(K)$.
Let $H^m(\Omega; X),
H^m(\T; X)$, resp.~$P_{\hspace{-.13em}k}(\T; X)$ denote the space of (piecewise) Sobolev functions resp.~polynomials with values
in $X=\mathbb R^2, \mathbb R^{2\times2}, \mathbb S\subset\mathbb R^{2\times2}$ (symmetric $2 \times 2$ matrices).

Let $\nabla_\pw\coloneqq D_\pw$, $D^2_\pw$, and $\div_\pw$ denote the piecewise gradient, Hessian, and divergence operators without explicit reference to the underlying triangulation $\T$.
Notice that 
 $ (H^2(\cT), a_\pw+j_h)$ becomes a Hilbert space~\cite[Sec.~4]{CCBGNN22} with
 the scalar product 
 $a_\pw+ j_h:H^2(\T)\times H^2(\T)\to \mathbb R$ defined by
\begin{align} 
a_\pw(v_\pw,w_\pw)&\coloneqq(D^2_\pw v_\pw, D^2_\pw w_\pw)_{L^2(\Omega)}\qquad\text{for any }v_\pw,w_\pw\in H^2(\T),\label{eqn:a_pw_def}\\
	j_h(v_{\pw},w_\pw)&\coloneqq
	\sum_{E \in \E}\left( \sum_{z \in {\mathcal V} (E)} \frac{[v_{\pw}]_E(z)}{h_E}\frac{[w_\pw]_E(z)}{h_E}
			 +
\fint_E  \jump{\frac{\partial v_{\pw}}{\partial\nu_E}}_E\mathrm ds\, \fint_E  \jump{\frac{\partial w_\pw}{\partial \nu_E}}_E\mathrm ds\right)\label{eqn:jh_defn}
\end{align}
with the jumps $\jump{v_{\pw}}_{{E}}(z) $
for $z \in \mathcal{V}(E)$ and $E \in \mathcal{E}$ defined as follows. The edge-patch $\omega(E):=\text{\rm int}(T_+\cup T_-)$ of an interior edge 
$E=\partial T_+\cap\partial T_-\in\E(\Omega)$  is the interior of the union 
$T_+\cup T_-$ of two neighboring triangles $T_+$ and $T_-$. Fix the orientation of the unit normal $\nu_E$ along $E$ and label $T_{\pm}$ such that $\nu_{T_+}|_E=\nu_E=-\nu_{T_-}|_E$ is the outer normal of $T_+$ along $E$. Let $\partial_s$ denote the tangential derivative along an edge $E$.
Then the jump resp.\ average read $[v_{\rm pw}]_E:=(v_{\pw} |_{T_+} -v_{\pw} |_{T_-})$ resp.\ $\langle v_{\pw} \rangle_E\coloneqq\half\left(v_{\pw}|_{T_+}+v_{\pw}|_{T_-}\right)$ on $E \in \E(\Omega)$. 
For a boundary edge $E \subset\partial T\cap \partial\Omega$ contained in the unique triangle $T \in \T$, $v_E= v_T|_E$, set $\omega(E) ={\rm int}(T)$ and $\jump{ v_{\pw}}_{E} = v_{\pw}|_E $ resp.\ $\langle v_{\pw}\rangle_E:=v_{\pw}|_E$. Let $\tau_E$ denote the unit tangent of fixed orientation along an edge $E\in\E$ and abbreviate $h_{\E},\nu_{\E},\tau_{\E},$ resp.\ $\jump{v_\pw}_{\E}$ as functions on the skeleton $\bigcup \E$ with $h_{\E}|_E\coloneqq h_E, \nu_{\E}|_E\coloneqq\nu_E,\tau_{\E}|_E\coloneqq\tau_E,$ resp.\ $\jump{v_\pw}_{\E}|_E:=\jump{v_\pw}_E$ for any $E\in\E$.
The piecewise integral mean operator $\Pi_{\E,0}$ reads $\Pi_{\E,0}(v)|_E\coloneqq\Pi_{E,0}(v)\coloneqq\fint_E v \;\ds$ for any $v\in L^2(E)$ and $E\in\E$.

\medskip \noindent 
The remaining parts of this paper apply the abstract results from Section \ref{sec:Abstract a posteriori error control} to fourth-order problems with the Sobolev spaces $V=X=Y:=H^2_0(\Omega)$
endowed with the energy norm $\trb{\bullet}\equiv(\trb{\bullet}_\pw)|_V\equiv(\|\bullet\|_h)|_V$ for the seminorm $\trb{\bullet}_\pw\coloneqq(a_\pw(\bullet,\bullet))^{1/2}  $ in $H^2(\T)$
and the discrete spaces $V_h=X_h=Y_h \subseteq P_2(\T)$ equipped with the induced norm $\|\bullet\|_h$ of the common superspace $\widehat V=\widehat X=\widehat Y:=H^2(\T)$ given as
\begin{align}
    \label{hnorm}
\|v_{\pw}\|_{h}^2 &:= \trinl v_{\pw} \trinr_{\pw}^2 + j_h(v_{\pw}, v_{\pw}) \qquad\text{for all } v_{\pw} \in H^2(\cT).
\end{align}
The subsequent analysis also requires the Morley finite element space
\[
\cM(\cT):=\left\{ v_\cM\in P_2(\cT) \; \vrule\;
\jump{v_\M}_E(z) = 0\text{ and }\int_{E}\jump{\frac{\partial v_\cM}{\partial \nu_E}}_E{\ds}=0\text{ for all }E\in\E\text{ and } z\in\mathcal V(E)
\right\}
\]
that lies in the kernel of $j_h$, i.e., $j_h(v_\M,\bullet)=0$ such that $\|v_\M\|_h=\trb{v_\M}_\pw$ for all $v_\M\in\M(\T)$, and comes with 
the Morley interpolation operator $I_\M$ that generalizes from $V$ to $H^2(\T)$ by averaging~\cite{ccnnlower2022}.
\begin{defn}[{Morley interpolation~\cite[Definition 3.5]{ccnnlower2022}}]   \label{def:morleyii}
	Given any $v_{\pw} \in H^2(\cT)$, define $I_{\rm M}v_{\pw} := v_{\rm M} \in {\rm M}(\cT)$ by its degrees of freedom as follows.  For any interior vertex $z \in \mathcal{V}(\Omega)$ with the set of attached triangles $\cT(z)$ of cardinality $|\cT(z)| \in \mathbb{N}$ and for any interior edge {{$E \in \mathcal{E}(\Omega)$}},
	\begin{align}
	v_{\rm M}(z) := |\cT(z)|^{-1} \sum_{K \in \cT(z)} (v_{\pw}\vert_{K})(z)\;\;\text{ and } \fint_{E} \dfrac{\partial v_{\rm M}}{\partial \nu _{\rm E}}\,\mathrm{d}s := \fint_{E} \left\langle \dfrac{\partial v_{\rm pw}}{\partial \nu_{E}} \right\rangle_E\,\mathrm{d}s.
	\end{align} 
	The remaining degrees of freedom at the vertices and edges on the boundary are set zero owing to the homogeneous boundary conditions.
\end{defn}
\noindent
An important property \cite[Eqn.~(3.5)]{ccnnlower2022} of the generalized Morley interpolation is the $a_\pw$-orthogonality
\begin{equation} \label{eqn:rightinverse}
    a_{\pw}(v_2, v-I_{\M} v)  =0 \qquad\text{for all } v \in V \text{ and all } v_2 \in P_2(\T).
\end{equation}
The point is that there exists a right-inverse $J:\M(\T)\to V$ of $I_\M$, that is, $I_{\M}J v_\M= v_\M$ for all $v_\M \in \M(\T)$, as in \cite{DG_Morley_Eigen,ccnn2021}, \cite[Lemma 3.7, Theorem 4.5]{ccnnlower2022} with
$%
    \|v_h-JI_\M v_h\|_h\lesssim\min_{v\in V}\|v-v_h\|_h
$
for all $v_h\in P_2(\T)$
such that $P=Q=S\coloneqq JI_\M: V_h \rightarrow V$ and $R \in \{{\rm id}, I_\M, JI_\M \}$ satisfy the assumptions \eqref{quasioptimalsmootherP}--\eqref{quasioptimalsmootherS} of Section \ref{sec:Abstract a posteriori error control}.
The last ingredient is a bounded transfer operator  $I_h: \M(\T) \rightarrow V_h$ 
that 
is either 
the identity $I_h\coloneqq{\rm  id}$ for the Morley, dG, and WOPSIP schemes with $\M(\T)\subset V_h$ or
$I_h\coloneqq I_{\rm C}$ for C$^0$IP with $V_h:=S^2_0(\T)\coloneqq P_2(\T) \cap H^1_0(\Omega)$ defined, for all $v_\M\in \M(\T)$, by
\begin{align} \label{eq:ic}
(I_{\rm C} v_\M)(z)=
\begin{cases}
	v_\M(z) &\text{for all } z \in {\mathcal V},\\
	\langle{v_\M}\rangle_E (z)  &\text{for } z= \text{mid}(E), \; E \in \E(\Omega), \\
	0 &\text{for } z= \text{mid}(E), \; E \in \E(\partial \Omega).
\end{cases}
\end{align}
The boundedness of $I_h$ implies $\|I_h\|\coloneqq\sup_{v_\M\in\M(\T)}\|I_h v_\M\|_h/\trb{v_\M}_\pw<\infty$.
\subsection{Explicit residual-based a posteriori estimator}\label{sub:explicit_aposteriori}
Any general source $F \in H^{-2}(\Omega)$ can be written with $L^2$ functions $f_0 \in L^2(\Omega), f_1 \in L^2(\Omega; {\mathbb R}^2), f_2 \in L^2(\Omega; {\mathbb S})$ \cite[Thm.~7.1]{CCBGNN22} as 
\begin{align} \label{eqn:newone}
F (\varphi) := \int_{\Omega} (f_0\; \varphi + f_1 \cdot \nabla \varphi + f_2 : D^2 \varphi) \dx \text{ for all } \varphi \in H^2_0(\Omega).
\end{align}
This definition extends to arguments $\varphi_\pw \in H^2(\T)$ by replacing $\varphi, \nabla \varphi, D^2 \varphi$ by their piecewise versions $\varphi_\pw, \nabla_\pw \varphi_\pw, D^2_\pw \varphi_\pw$. In the applications below, the approximation $v_h \in V_h \subset P_2(\T)$ is fixed 
 and $\Gamma_{\pw}(v_h,v_h,\bullet) \in H^{-2}(\Omega) $ has a structure as in \eqref{eqn:newone}, namely
\begin{align}\label{eqn:newtwo}
\Gamma_{\pw}(v_h,v_h, \varphi) = \int_{\Omega} (\Gamma_0\; \varphi + \Gamma_1 \cdot \nabla \varphi 
+ \Gamma_2: D^2 \varphi) \dx  \text{ for all } \varphi \in H^2_0(\Omega)
\end{align} 
for piecewise polynomials $\Gamma_0 \in P_k(\T)$, $\Gamma_1 \in P_k(\cT; {\mathbb R}^2), $
 $\Gamma_2 \in P_k(\cT; {\mathbb S})$ of degree at most $k\in\mathbb N_0$.
 The Lebesgue functions in \eqref{eqn:newone} and the polynomial degree $k\in\mathbb N_0$ give rise to oscillations 
 \begin{align} \label{eqn:newfour}
 {\rm osc}_k(F,\T) &:= \| h_\T^2(f_0-\Pi_k f_0)\|_{L^2(\Omega)} + \| h_\T(f_1-\Pi_k f_1)\|_{L^2(\Omega)} 
 + \| f_2-\Pi_k f_2\|_{L^2(\Omega)}.
 \end{align}

\begin{example}[Navier-Stokes] \label{ex:3.1} The trilinear form $\Gamma_\pw$ for the Navier-Stokes equations in Section \ref{sec:nse} below is given by \eqref{eqn:newtwo} for $k=1$ with $\Gamma_1= \Delta_{\pw} v_h \: {\rm Curl}_{\pw} v_h$ and $\Gamma_0=0$, $\Gamma_2=0$.
\end{example}
\begin{example}[von K\'{a}rm\'{a}n]\label{ex:3.2}
For the von K\'{a}rm\'{a}n equiations in Section~\ref{sec:vke} and vector-valued approximation $\mathbf{v_h}\equiv\big(\vha,\vhb\big)\in V_h\times V_h$, 
 choose $\Gamma_0=[\vha,\vhb]$  resp.\  $\Gamma_0=-1/2[\vha,\vha]$ and $\Gamma_1= 0$, $\Gamma_2=0$ with the von K\'{a}rm\'{a}n  bracket 
$[\bullet,\bullet]$
defined in Subsection \ref{vke:model_intro} and $k=0$. 
\end{example}

 \medskip
\noindent Recall $\tau_{\E},\nu_{\E},[\bullet]_{\E}$, and the piecewise integral mean operator $\Pi_{\E,0}$ from Subsection \ref{sub:triangulation, interpolation, and smoother}. The error estimators $ \mu_1(\T) + \mu_2(\T) + \mu_3(\T)=: \mu(\T) $ and $\eta(\T)$ are defined in terms of 
\begin{equation}\label{eqn:lambdas}
\begin{array}{ll}
{\Lambda}_0:= \Pi_k f_0 - \Gamma_0 \in P_k(\T), \\
\Lambda_1:= \Pi_k f_1 - \Gamma_1 \in  P_k(\cT; {\mathbb R}^2),  \\  
\Lambda_2 := \Pi_k f_2 - D^2_{\pw} v_h - \Gamma_2 \in P_k(\cT; {\mathbb S})
\end{array}
\end{equation}
with $\vartheta=1$ for $I_h={\rm id}$ resp. $\vartheta=0$ for $I_h=I_{\rm C}$ by 
 \begin{align*}
 \mu_1(\T) &:= \big\| h_{\T}^2 (\Lambda_0 - {\rm div}_{\pw} \Lambda_1 +{\rm div}^2_{\pw} 
 \Lambda_2) \big\|_{L^2(\Omega)}, \\
  \mu_2(\T) &:= \big\| h_{\E}^{3/2}\jump{\Lambda_1 - {\rm div}_{\pw} \: 
 \Lambda_2 - \partial_s (\Lambda_2 \tau_{\E})}_{\E} \cdot \nu_{\E}  \big\|_{L^2(\E(\Omega))},  \\
  \mu_3(\T) &:= 
  \big\|h_{\E}^{1/2} (1-\vartheta \: \Pi_{{\E},0}) \jump{\Lambda_2 \nu_{\E}}_{\E} \cdot  \nu_{\E} \big\|_{L^2(\E(\Omega))}, 
\\
\eta(\T)&:=\| v_h-{JI_\M}  v_h\|_{h}.
 \end{align*}
The following theorem controls the intermediate residual $\rho(M)$ from \eqref{a} with $M=\|I_h\|$ from the abstract reliability estimate \eqref{eqn:apost} by the explicit a posteriori error estimators $\mu(\T)$ and $\eta(\T)$. Define
\begin{align}\label{eqn:varrho}
    \rho\coloneqq
\sup_{\stackrel{y \in Y, \|y \|_Y = 1}{z:=y- J I_\M I_hI_\M y} } 
\big( F(z)-a( J I_\M v_h,z) -
\Gamma_{\pw}(Rv_h, Rv_h, z)  \big).
\end{align}

 \begin{theorem}[reliability and efficiency up to data oscillations] \label{thm:releff}
 Under the present notation $V=X=Y\equiv H^2_0(\Omega)$, the choices $P=Q=S=JI_\M$, and
 $R=\{{\rm id}, I_\M, J I_\M \}$, it holds 
 \begin{align*}
& (a)   \; \rho(\|I_h\|) \leq\rho\lesssim  \eta(\T) + \mu(\T) +{\rm osc}_k(F,\T), \\
&
 (b) \;  \eta(\T) + \mu(\T)\lesssim  \rho(0)+ {\rm osc}_k(F,\T), 
 \\
 & (c) \; \eta(\T) \approx \big\|h_{\E}^{1/2} \jump{D^2_\pw v_h}_{\E}\tau_{\E}\big\|_{L^2(\E)} + j_h(v_h,
		v_h)^{1/2}.
  \end{align*}
\end{theorem}
\noindent The estimate Theorem~\ref{thm:releff}.c is well-known from \cite[Theorem 5.6]{CCBGNN22}, \cite{DG_Morley_Eigen}. The remaining parts of this section therefore focus on the proofs of Theorem~\ref{thm:releff}.a and b.
\begin{rem}[role of \eqref{eqn:QO}] The a priori results in \cite{CCNNGRDS22} establish \eqref{eqn:QO} for  $S=Q=JI_\M$ in \eqref{eqn:discrete} and this leads to efficiency. Theorem \ref{thm:apost} is fairly general and the reliability estimate allows for $S \neq Q$; but then \eqref{eqn:QO} involves an additional additive term $\Gamma(Ru_h, Ru_h, (S-Q)y_h) = O(h_{\rm max}^{\alpha})$ \cite[Theorem 5.1]{CCNNGRDS22}. This extra term behaves like a given $L^2$ function (in terms of piecewise derivatives of $Ru_h$)  times the mesh-sizes up to some power  $\alpha \ge 0$. The application to Navier-Stokes leads to $\alpha=1$ and this is of the correct asymptotic rate (or even better), while the application to \vket ~even allows for $\alpha=2$ {\cite{CCNNGRDS22}}. However, this extra term is {\it not} a residual term (in general) and efficiency is left open as in \cite{CCGMNN_semilinear,kim_morley_2021}. 
The new schemes from \cite{CCNNGRDS22} with $P=Q=S=JI_\M$ in \eqref{eqn:discrete}  enable an efficient and reliable a posteriori error control in this paper for general sources.  
\end{rem}

\begin{comment}
{\begin{theorem} \cite{CCBGNN22} \label{thm:fromold}There exist positive constants $C_{\rm rel}^{\ell}, C_{\rm eff}^{{\ell}}$, that exclusively depend o the shape regularity of $\T$ and on the polynomial degree $k \in {\mathbb N}_0$ such that 
 $C_{\rm rel}^{\ell} \trinl \Lambda \circ (1-J_hI_h I_\M) \trinr_* \le \mu(\T) \le C_{\rm eff}^{{\ell}} \trinl \Lambda \trinr_*$.
 \end{theorem}
 \begin{lemma} \cite{CCBGNN22} \label{lem:newone} Any $v \in H^2_0(\Omega)$ satisfies $$ \sum_{j=0}^2 \|h_\T^{j-2} (v- JI_\M I_h I_\M v) \| \lesssim \mu(\T) \lesssim \trinl  \Lambda \trinr_*.$$
 \end{lemma}
 \begin{lemma} \cite{CCBGNN22} \label{lem:newtwo} $\trinl v_h- JI_\M v_h \trinr_{\pw} \approx \eta(\T) \lesssim \|v_h-Jv_h\|_h \lesssim \trinl u-v_h\trinr_{\pw}$.
 \end{lemma} }
 \footnote{to be checked and quoted properly}
\end{comment}
\subsection{Proof of Theorem~\ref{thm:releff}.a.} \label{appl:rel}
The definition of $\rho(M)$ for $M\coloneqq\|I_h\|$ in \eqref{a} implies $\rho(M)\leq\rho$.
Indeed, for any $v\in V$ with $\trb{v}=1$ and $v_h:= I_h I_\M v \in V_h$, it follows that $\|v_h\|_h \le \|I_h\| \trinl  I_\M v \trinr_{\pw} \le M$ from 
$ \trinl I_\M v\trinr_{\pw} \le \trinl v \trinr = 1$ and so $v_h\equiv I_hI_\M v$ is admissible (the last inequality is a consequence of the Pythogoras theorem and the orthogonality \eqref{eqn:rightinverse}). %
The reflexivity of $V \equiv H^2_0(\Omega)$ leads to $w \in V$ with $\trinl w \trinr = 1$ and
\begin{equation}\label{eqn:star1}
    \rho(\|I_h\|)  \leq \rho \le   F(z)-a(Pv_h,z)-\Gamma_{\pw}(Rv_h,Rv_h,z) 
\end{equation}
for $z:=w-JI_\M w_h$. 
Recall $\Lambda_0,\Lambda_1, \Lambda_2$ from \eqref{eqn:lambdas} and define $\Lambda \in H^{-2}(\Omega)$ by 
 \begin{equation} \label{eqn:newthree}
 \Lambda(\varphi)= \int_{\Omega} (\Lambda_0 \varphi + \Lambda_1 \cdot \nabla \varphi + \Lambda_2: D^2 \varphi) \dx \text{ for all } \varphi \in H^2_0(\Omega).
 \end{equation}
Observe carefully that the definition of $F$ in \eqref{eqn:newone} and $\Lambda$ in \eqref{eqn:newthree}  lead to 
\begin{align}\label{eqn:startwo}
F-a(Pv_h,\bullet)-\Gamma_{\pw}(Rv_h,Rv_h,\bullet) 
& = F- \Pi_k F + \Lambda + a_{\pw}(v_h-Pv_h,\bullet) \nonumber \\
& \quad  +  \Gamma_{\pw}(v_h,v_h,\bullet)-\Gamma_{\pw}(Rv_h,Rv_h,\bullet).
\end{align}
Here and throughout, $\Pi_k F \in H^2(\T)^*$ is defined by
\begin{align}\label{eqn:Pi_F_def}
    \Pi_k F (v_\pw) := 
\int_{\Omega} \big( v_\pw \Pi_k f_0 +
\nabla_{\pw}v_\pw \cdot \Pi_k f_1+ D^2_{\pw} v_{\pw}:
\Pi_k f_2 \big) \dx \quad\text{for } v_{\pw} \in H^2(\T).
\end{align}
The six terms on the right-hand side in \eqref{eqn:startwo} are controlled as follows. Since $z=(1-J I_\M I_h I_\M )w$ vanishes at the vertices for all $w \in V$, the stability result \cite[Lemma 5.1]{CCBGNN22} 
 \begin{align} \label{eqn:C_stability}
 \sum_{m=0}^2 |h_\T^{m-2} (1- JI_\M I_h I_\M )w |_{H^m(\T)} & \le C_{\rm stab} \trinl w  \trinr \text{ for all } w \in H^2_0(\Omega)
 \end{align}
controls the  data oscillation term  
\cite[Lemma 7.2]{CCBGNN22} by  
\begin{align} \label{eqn:fbound}
   (F- \Pi_k F)(z) &= \int_{\Omega} \Big( 
   (f_0- \Pi_k f_0) z + 
    (f_1- \Pi_k f_1) \cdot \nabla z +
   (f_2- \Pi_k f_2): D^2 z \Big) \dx \nonumber \\
   & 
    \le C_{\rm stab} {\rm osc}_k(F,\T) \trinl z \trinr. 
\end{align}

\medskip \noindent A Cauchy inequality, the boundedness of $\Gamma_{\rm pw}$, $\trinl v_h-Pv_h \trinr_{\pw} \leq \|v_h-Pv_h\|_h$ from \eqref{hnorm}, and the definition of $\eta(\T) $ %
 reveal 
    \begin{align}
    a_\pw(v_h-Pv_h,z) & \le  \trinl v_h-Pv_h \trinr_{\pw} 
    \trinl z \trinr \le \eta(\T) \trinl z \trinr, \\
 \Gamma_{\pw}(v_h,v_h,z)-\Gamma_{\pw}(Rv_h,Rv_h,z)  & \le \|\Gamma_\pw\| (1+\|R\|) \|v_h\|_{h}
\|v_h-Rv_h\|_{h} \trinl z \trinr \nonumber \\
 &\le \|\Gamma_\pw\| (1+\|R\|)\lamr\big(\trb{u} +  \beta/ ((1+\Lambda_{\rm P}) \|\Gamma\|)\big)\; \eta(\T)  \trinl z \trinr\label{eqn:G_bound2a}
\end{align}
with the arguments from the analysis of $S_2$ in Subsection~\ref{S2}, 
$\trinl v_h-Rv_h \trinr_{\pw} \leq \|v_h-Rv_h\|_h$ from \eqref{hnorm}, \eqref{quasioptimalsmootherP}, \eqref{quasioptimalsmootherR}, and \eqref{eqn:vhbound}  in the last step.
It remains to control $\Lambda$,   where we employ results from the linear situation. 
The appendix of the preliminary work on linear problems \cite{CCBGNN22} provides the estimate 
\begin{align} \label{eqn:starfour}
\trinl \Lambda \circ (1-JI_\M I_h I_\M) \trinr_* \le 
C_{\rm rel}^{\rm lin} \: \mu(\T).
 \end{align}
 Since $z \equiv w- JI_\M I_h I_\M w$ acts as a test function in \eqref{eqn:starfour} and $\trb{w}=1$, we infer 
 \begin{align} \label{eqn:lam}
 \Lambda(z) \leq\trinl \Lambda \circ (1-JI_\M I_h I_\M) \trinr_* \trb{w}\le C_{\rm rel}^{\rm lin}\: \mu(\T).
 \end{align}
   Since $\trinl z \trinr \le C_{\rm stab}$ from \eqref{eqn:C_stability},
   the reliability  $\rho(\|I_h\|) \leq\rho\lesssim  \eta(\T) + \mu(\T) +{\rm osc}_k(F,\T)$ follows from \eqref{eqn:startwo},\eqref{eqn:fbound}--\eqref{eqn:G_bound2a}, and \eqref{eqn:lam} in \eqref{eqn:star1}. 
   The above constants $C_{\rm stab}$ and $C_{\rm rel}^{\rm lin} $ exclusively depend on the shape regularity of the triangulation $\T$ and the polynomial degree $k$ of the $\Lambda_0,\Lambda_1,\Lambda_2$ in \eqref{eqn:newtwo}. 
   \qed
\subsection{Proof of Theorem~\ref{thm:releff}.b} 
The efficiency $\eta(\T) \equiv \|v_h-Pv_h\|_h \le \Lambda_{\rm P }\| u-v_h \|_{h} $ follows from \eqref{quasioptimalsmootherP} for $P\equiv JI_\M$ as in \cite[Thm~5.6]{CCBGNN22}.%

\medskip \noindent The efficiency of $\mu(\T) \le  C_{\rm eff}^{\rm lin}
\trinl \Lambda \trinr_*$ is established in \cite[Theorem A.1]{CCBGNN22}  with $C_{\rm eff}^{\rm lin}$ that exclusively depends on the shape-regularity of $\T$. The definitions of $\Lambda= \Pi_k F - a(v_h,\bullet ) - \Gamma_{\pw} (v_h,v_h,\bullet)$ and $\rho(0)$ therefore lead to 
\begin{align} \label{eqn:mubound}
\mu&(\T)  \le C_{\rm eff}^{\rm lin} \trinl \Lambda \trinr_*  \le  \trinl \Pi_k F-a(v_h,\bullet)-\Gamma_{\pw}(v_h,v_h,\bullet) \trinr_* \nonumber \\
&\le  \rho(0) + \trinl F- \Pi_k F\trinr_* + \trinl a_{\pw} (v_h - Pv_h,\bullet) \trinr_* 
  + \trinl \Gamma_{\pw}(v_h,v_h,\bullet)-\Gamma_{\pw}(Rv_h,Rv_h,\bullet) \trinr_*.
\end{align}
A Cauchy inequality and the boundedness and $\Gamma_{\pw}$ provide, for any $z\in V$, that
  \begin{align} \label{bounda}
    a_\pw(v_h-Pv_h,z) & \le \trb{ v_h-Pv_h }_\pw\trb{z}
    \le \Lambda_{\rm P} \|u-v_h\|_{h}\trinl z \trinr,  \\ 
 \Gamma_{\pw}(v_h,v_h,z)\!-\!\Gamma_{\pw}(Rv_h,Rv_h,z)  & \le \|\Gamma_\pw\| (1\!+\!\|R\|) \|v_h\|_{h}
\|v_h-Rv_h\|_{h} \trinl z \trinr \nonumber \\
 &\le \|\Gamma_\pw\| (1\!+\!\|R\|)\lamr\big(\trb{u} \!+ \! \beta/ ((1\!+\!\Lambda_{\rm P}) \|\Gamma\|)\big) 
\|u\!-\!v_h\|_{h} \trinl z \trinr \label{eqn:bound_gamma}
\end{align}
with $\trinl v_h-Pv_h \trinr_{\pw} \leq \|v_h-Pv_h\|_h$ from \eqref{hnorm}, \eqref{quasioptimalsmootherP} for  \eqref{bounda} and $\trinl v_h-Rv_h \trinr_{\pw} \leq \|v_h-Rv_h\|_h$ from \eqref{hnorm}, \eqref{quasioptimalsmootherR},  and \eqref{eqn:vhbound} in the last step. 
The combination of \eqref{eqn:mubound}--\eqref{eqn:bound_gamma} and \eqref{eqn:fbound}
lead to 
$\mu(\T) \lesssim \rho(0) + {\rm osc}_k(F,\T) + \|u-v_h \|_h$ and conclude the proof.
\qed

\section{Application to Navier-Stokes equations} \label{sec:nse}
\subsection{Stream-function vorticity formulation of Navier-Stokes equations} \label{sec:prob}
The stream-function vorticity formulation of the incompressible 2D Navier--Stokes problem for a given load $F\in H^{-2}(\Omega)$ in a polygonal domain $\O \subset \R^2$ seeks $u \in V:=H^2_0(\O) (\equiv X\equiv Y)$  such that%
\begin{equation} \label{NS_eqa}
 \Delta^2u + \frac{\partial}{\partial x}\Big((-\Delta u)\frac{\partial u}{\partial y}\Big)- \frac{\partial}{\partial y}\Big((-\Delta u)\frac{\partial u}{\partial x}\Big)= F \quad\text{ in } \O.
\end{equation}
(The bi-Laplacian $\Delta^2$ reads $\Delta^2\phi:=\phi_{xxxx}+\phi_{yyyy}+2\phi_{xxyy}$.)  The analysis of extreme viscosities lies beyond the scope of this article, and so the viscosity of the
 bi-Laplacian in \eqref{NS_eqa} is set one. 
 Recall the semi-scalar product $a_{\rm pw}:H^2(\T)\times H^2(\T)\to\mathbb R$ and the induced piecewise $H^2$ seminorm $\trinl \bullet \trinr_{\rm pw} $ that is a norm \cite{ccnn2021} on $V+\M(\cT) \subset H^2(\T)$ from Subsection \ref{sub:triangulation, interpolation, and smoother}.
 Define the bounded trilinear form $\Gamma_{\pw}(\bullet,\bullet,\bullet)$ by 
\begin{align} \label{deftns:bilinear}
\Gamma_{\pw}(\widehat{\eta}, \widehat{\chi},\widehat{\phi}) :=\sum_{T\in\T}\int_{T}^{} \Delta \widehat{\eta}\,\Big(\frac{\partial \widehat{\chi}}{\partial y}\frac{\partial \widehat{\phi}}{\partial x}-\frac{\partial \widehat{\chi}}{\partial x}\frac{\partial \widehat{\phi}}{\partial y}\Big) \dx\qquad\text{for all }\widehat{\eta},\widehat{\chi},\widehat{\phi} \in   H^2(\T).
\end{align} 
Given $a:=a_{\pw}|_{V \times V}$ and $\Gamma:=\Gamma_{\pw}|_{V\times V \times V}$, 
the weak formulation of  \eqref{NS_eqa} seeks $u \in V$ such that
\begin{align} \label{NS_weak}
a(u,w) + \Gamma(u,u,w) = F(w)\fl w \in V.
\end{align}
Given any $F\in V^*\equiv H^{-2}(\Omega)$, there exist solutions to \eqref{NS_weak}, which are possibly singular but carry elliptic regularity. 
In the case of small loads ($\trb{F}_*\|\Gamma\|<1$), the weak solution is unique and a regular root, cf.~\cite[Chap.~IV.§2--3]{GR:FiniteElementMethods1986} and \cite{NS_FE79,Tem:NavierStokesEquations1995} for proofs.
The a posteriori error analysis below concerns some approximation $v_h \in  V_h$  to a regular root  ${u} \in V$ of the continuous problem~\eqref{NS_weak}.}
\subsection{Five quadratic discretizations}\label{sec:discrete_schemes}
This subsection presents the {Morley},  two variants of {dG},  ${ C^0}$IP, and  WOPSIP discretizations for \eqref{NS_weak}.
The discrete space $V_h\equiv X_h\equiv Y_h$ becomes $V_h\coloneqq\M(\T)$  for Morley, $V_h\coloneqq P_2(\T)$ for dG and WOPSIP schemes, and $V_h\coloneqq S^2_0(\cT)\coloneqq P_2(\cT) \cap H^1_0(\Omega)$ for the  $C^0$IP scheme. For all  $v_{\rm pw}, w_{\rm pw} \in H^2(\cT)$ and  parameters $\sigma_1,\sigma_2, \sigma_{\ip}>0$ sufficiently large (but fixed in applications) to guarantee the stability of $a_h$ below, the method-dependent penalty forms $c_{\dg}, c_{\rm P}$, and $c_{\rm IP}$ read
\begin{align}\label{eqn:cdG_def}
    c_{\rm dG}(v_{\pw}, w_{\pw})&\coloneqq \sum_{E\in\cE}\bigg(\frac{\sigma_1}{h_E^3}\int_E \jump{v_{\rm pw}}_E \jump{w_{\rm pw}}_E\ds+ \frac{\sigma_2}{h_E}\int_E\jump{ \frac{\partial v_{\rm pw}}{\partial  \nu_E}}_E  \jump{\frac{\partial w_{\rm pw}}{\partial  \nu_E}}_E \ds\bigg),\\
    c_{\rm P}(v_{\pw}, w_{\pw})&\coloneqq \sum_{E \in \mathcal{E}} 
h_{E}^{-4} \bigg(\sum_{z \in \mathcal{V}(E)} \Big(\jump{v_{\rm pw}}_E   
\jump{w_{\rm pw}}_E\Big)(z)+\int_E\jump{ \frac{\partial v_{\rm pw}}{\partial  \nu_E}}_E  \ds \int_E\jump{\frac{\partial w_{\rm pw}}{\partial  \nu_E}}_E \ds\bigg),\label{eqn:c_P_def}
\end{align}
and $c_{\rm IP}\coloneqq c_{\rm dG}|_{(V+S^2_0(\T))\times (V+S^2_0(\T))}$ with $\sigma_{\rm IP}\coloneqq \sigma_{2}$. Define the discrete bilinear forms 
\[a_h:= a_{\rm pw} + {\mathsf b}_h+ 
{\mathsf c}_h:\left(V + V_h\right)\times \left(V + V_h\right)\to \mathbb R, \] with $a_{\pw}$ from \eqref{deftns:bilinear} for the Morley, dG I, ${ C^0}$IP, and  WOPSIP discretizations and $a_{\pw}$ replaced by $(\Delta_{\pw} \bullet, \Delta_{\pw} \bullet)_{L^2(\Omega)}$
for the dG II scheme, and  ${\mathsf b}_h$ and ${\mathsf c}_h $ from Table \ref{tab:spaces} for some $-1\leq\theta\leq 1$.
The method-dependent norms induced by $a_\pw+c_h$ for dG I, WOPSIP, and $C^0$IP {(resp.  
$(\Delta_\pw\bullet,\Delta_\pw\bullet)_{L^2(\O)} +c_h$ for dG II as in \cite{SM:HpversionInteriorPenalty2007})
are, except for WOPSIP, equivalent to the universal norm $\|\bullet \|_{h}$ from 
\eqref{hnorm}. Notice that $\|\bullet\|_h=\trb{\bullet}_\pw$ in $V+\M(\T)$ follows from \eqref{hnorm}.
\begin{lemma}[{Equivalence of norms \cite[Thm.~4.1]{ccnngal}}] \label{lem:equivalence}
	It holds 
 $\| \bullet \|_{h} \approx \|\bullet \|_{\rm dG}\equiv \big(\trinl \bullet \trinr_{\rm pw}^2 +c_{\dg}(\bullet,\bullet)\big)^{1/2}$ on $V + P_2(\cT)$ and $\| \bullet \|_{h} \approx \|\bullet \|_{\rm IP}\equiv \big(\trinl \bullet \trinr_{\rm pw}^2 +c_{\rm IP}(\bullet,\bullet)\big)^{1/2}$ on $V + S^2_0(\cT)$. \qed
\end{lemma}  
In contrast to this, the WOPSIP norm $\|\bullet \|_{\rm P} \equiv \big(\trinl \bullet \trinr_{\rm pw}^2 +c_{\rm P}(\bullet,\bullet)\big)^{1/2}$ involves smaller powers of the mesh-size and is (strictly) stronger than $\|\bullet\|_h$, i.e., $h_T\leq h_{\rm max}$ implies
\begin{align}\label{eqn:P_norm}
    j_h\leq h_{\rm max}^2 c_{\rm P}\qquad\text{and}\qquad \|\bullet\|_h\leq (1+h_{\rm max}^2)^{1/2}\|\bullet\|_{\rm P}.
\end{align}
\noindent  The applications in this paper consider the choice $P\equiv Q\equiv S\equiv JI_\M$  %
that allows the first reliable and efficient a posteriori error estimate for the stream-function vorticity formulation of the Navier-Stokes equations. 
\begin{table}[]
\centering {\footnotesize
\begin{tabular}{|c|c|c|c|c|c|}
\hline
Scheme & Morley & \multicolumn{1}{c|} {dG I} &  $C^0$IP & WOPSIP & dG II \\ \hline

$\mathcal{J}(\bullet ,\bullet)  $ & -- & \multicolumn{2}{c|}{{${{{}\displaystyle \sum_{E\in\cE}\int_E \langle D^2v_2\;\nu_E\rangle_E} \cdot \jump{\nabla w_2}_E\ds}$}}  & -- &${{{}\displaystyle \sum_{E\in\cE}\int_E \jump{\frac{\partial v_2}{\partial \nu_E}}_E\langle \Delta_{\pw}w_2\rangle_E} }\,\ds$   \\\hline

${\mathsf b}_h(\bullet,\bullet)$ & 0 &
\multicolumn{2}{c|}{\begin{minipage}{3.6cm}
\centering
\vspace{3pt}
$ -\theta \mathcal{J}(v_2,w_2) {- \mathcal{J}(w_2,v_2)} $
\end{minipage}}
& 0 & $-\theta \mathcal{J}(v_2,w_2) {- \mathcal{J}(w_2,v_2)}$  \\ \hline
${\mathsf c}_h(\bullet,\bullet)$ & 0 &
$ c_\dg(\bullet,\bullet)$ & $ c_{\rm IP}(\bullet,\bullet)$ &$ c_{\rm P}(\bullet,\bullet)$ &$ c_\dg(\bullet,\bullet)$ \\\hline
$I_{h}$ & ${\rm id}$ & \multicolumn{1}{c|}{${\rm id}$} & $I_{\rm C}$  from \eqref{eq:ic} & ${\rm id}$ & ${\rm id}$ \\ \hline  
\end{tabular}
\caption{Bilinear forms with $\theta\in[-1,1]$ and operator $I_h$ in Section \ref{sec:nse}. }
\label{tab:spaces}}
\end{table}
Recall $\tau_{\E},\nu_{\E},[\bullet]_{\E}$, and the piecewise integral mean operator $\Pi_{\E,0}$ from Subsection~\ref{sub:triangulation, interpolation, and smoother} and abbreviate $\mathrm{Curl}\coloneqq(\partial/\partial y; - \partial/\partial x)$.
Given $\vartheta=1$ for $I_h={\rm id}$ resp. 
$\vartheta=0$ for $I_h=I_{\rm C}$, the local error estimators on $T\in\T$,
 \begin{align*}
 \sigma^2(T) & := |T|^2 \|\Pi_k f_0-  {\rm div}_{\pw} (\Pi_k f_1) +{\rm div}^2_{\pw}  (\Pi_k f_2) \|^2_{L^2(T)} \\
& \quad + |T|^{3/2} \|\jump{\Pi_k f_1 -\Delta v_h \: {\rm Curl} v_h -  {\rm div} (\Pi_k f_2) -
\partial_s ( \Pi_k f_2) \tau_{\E}}_{\E}  \cdot \nu_{\E} \|^2_{L^2(\partial T \setminus \partial \Omega)}\\
& \quad + |T|^{1/2}  { \| (1-\vartheta\Pi_{\E,0}) \jump{(\Pi_k f_2-  D^2 v_h)\nu_{\E}}_{\E} \cdot  \nu_{\E} \|_{L^2(\partial T \setminus \partial \Omega)}^2} 
\\
&\quad+|T|^{1/2} \big\|[D^2v_h]_{\E}\tau_{\E}\big\|_{L^2(\partial T)}^2+\sum_{E \in \mathcal{E}(T)} \Big(\big| \Pi_{E,0} \jump{\partial_\nu v_h}_E \big|^2 + |T|^{-1}\hspace{-0.6em} \sum_{z \in \mathcal{V}(E)} \big|\jump{v_h}_E(z)\big|^2 \Big)
     \end{align*}
define the a posteriori error estimator $\sigma(\T) := \sqrt{\sum_{T \in \T} \sigma^2(T)}$ by the $\ell^2$ sum convention.
\begin{rem}[Classical case]
    The a posteriori error analysis in \cite{CCGMNN_semilinear,kim_morley_2021} for the classical situation $R=S=Q={\rm id}$ and $F\equiv f\in L^2(\Omega)$ suffers from unclear efficiency.
\end{rem}
\begin{theorem}[A posteriori error control] \label{thm:apost_ns_final}
Given a regular root $u \in V$ to \eqref{NS_weak} with $F\in H^{-2}(\Omega)$ and $k\in {\mathbb N}$, there exist $\varepsilon, \delta,\varrho 
>0$ such that the following holds for any $\displaystyle\cT\in\bT(\delta)$ and $P\equiv Q\equiv S\equiv JI_\M$.
$(a)$ There exists a unique discrete solution $u_h\in V_h$ to \eqref{eqn:discrete} for the Morley, dG I \& II, and $C^0$IP scheme in Table~\ref{tab:spaces} with $\|  u-u_h\|_{h} \le\varepsilon$ and any $v_h \in V_h$ with $\|u_h-v_h\|_h \le \varrho$ satisfies
$$\| u-v_h \|_{h}  + {\rm osc}_k(F,\T) \approx \sigma(\T)  + \|u_h-v_h\|_{h} + {\rm osc}_k(F,\T).$$
$(b)$ For the WOPSIP scheme in Table~\ref{tab:spaces}, there exists a unique discrete solution $u_h\in V_h$ to \eqref{eqn:discrete} with $\| u-u_h\|_{h} \le\varepsilon$ and any $v_h \in V_h$ with $\|u_h-v_h\|_h \le \varrho$ satisfies
\begin{align}\label{eqn:WOPSIP_NVS_a}
    \| u-v_h \|_{h}  +\|u_h-v_h\|_{h} +  {\rm osc}_k(F,\T) &\approx \sigma(\T)  + \|u_h-v_h\|_{h} + {\rm osc}_k(F,\T),\\
    \| u-v_h \|_{\rm P}+\|u_h-v_h\|_{h}  +{\rm osc}_k(F,\T) &\approx \sigma(\T)  +c_{\rm P}(v_h,v_h)^{1/2}+ \|u_h-v_h\|_{h} + {\rm osc}_k(F,\T).\label{eqn:WOPSIP_NVS_b}
\end{align}
\end{theorem}
\noindent The Morley, dG I \& II, and $C^0$IP schemes satisfy the discrete consistency \eqref{eqn:H} and quasi-optimality \eqref{eqn:QO} so that the proof of Theorem \ref{thm:apost_ns_final} is already prepared in Theorems \ref{thm:apost} and \ref{thm:releff}.
The proof for the WOPSIP method (without \eqref{eqn:H} and \eqref{eqn:QO}) requires modifications in Subsection~\ref{sub:WOPSIP_NSV} below.
\begin{proof}[Proof of Theorem \ref{thm:apost_ns_final}.a]
The a priori analysis \cite[Thm.~8.1]{CCNNGRDS22} verifies the quasi-optimality \eqref{eqn:QO} and provides universal constants $\varepsilon_0,\delta_0>0$ that guarantee, for any $\T\in\TT{\delta_0}$, the unique existence of a discrete solution $u_h\in V_h$ to \eqref{eqn:discrete} with $\|u-u_h\|_h\leq \varepsilon_0$. 
A density argument for $\varepsilon\coloneqq\min\{\varepsilon_0,\beta/\left(3(1+\LP)\|\Gamma\|\right)\} $ leads to $\delta\leq\delta_0$ such that $\|u-u_h\|_h<\varepsilon$ for any $\T\in\TT\delta$.
This reveals \eqref{eqn:QO} and {\bf(L1)--(L3)} for $\delta,\varepsilon,\kappa\coloneqq2/3$, and $\varrho\coloneqq\beta/\left(3(1+\LP)\|\Gamma\|\right)$.
The abstract a posteriori error control from 
Theorem~\ref{thm:apost} applies with the abstract a posteriori error control of $\|u-v_h\|_h$ in terms of $\rho(\|I_h\|), \|(1-JI_\M)v_h\|_h$, and the efficient algebraic error $\|u_h-v_h\|_h$.
Section \ref{subsection:new} applied to $\Gamma_\pw$ from \eqref{deftns:bilinear} leads to 
${\Lambda}_0:= \Pi_k f_0, \; \Lambda_1:= \Pi_k f_1 -  \Delta_{\pw} v_h \: {\rm Curl}_{\pw} v_h, \; 
\Lambda_2 := \Pi_k f_2 - D^2_{\pw} v_h$ in \eqref{eqn:lambdas}.
Theorem~\ref{thm:releff}.a controls $\rho(\|I_h\|)+\|(1-JI_\M)v_h\|_h$ with the a posteriori term $\mu(\T) +\eta(\T)$, that is efficient by Theorem~\ref{thm:releff}.b and Theorem~\ref{thm:apost}.b, plus data oscillations ${\rm osc}_k(F,\T)$. The equivalence $\mu(\T) +\eta(\T)\approx \sigma(\T)$ follows with $h_E\approx h_T\approx |T|^{1/2}$ from shape-regularity and Theorem~\ref{thm:releff}.c. This concludes the proof.
\end{proof}
\subsection{Modifications for WOPSIP}\label{sub:WOPSIP_NSV}

There are two reasons why the weakly over-penalized symmetric
interior penalty (WOPSIP) scheme from \cite{BrenGudiSung10} requires little modifications in the above analysis.
The first is the failure of \eqref{eqn:H} and (somehow related) the failure of \eqref{eqn:QO} in the stated  form. The second is that the natural
WOPSIP norm $\|\bullet\|_{\rm P}$ from \eqref{eqn:P_norm} is very strong and \eqref{eqn:WOPSIP_NVS_a}
states the error estimate in the (partly) weaker norm $\|\bullet\|_h$.

\medskip 
\noindent 
The starting point for the analysis in this subsection is a modified version of \eqref{eqn:H} already applied in the analysis of linear problems \cite[Thm.~6.9]{CCBGNN22} that follows from $j_h(v_{\pw}, w_{\M})=0=c_{\rm P}(v_{\pw}, w_{\M})$ for all $(v_{\pw},v_{\M})\in H^2(\T)\times \M(\T)$.
Recall that $\|J\|$ abbreviates the operator norm of
$J: \M(\T)\to V$ when $\M(\T)$ and $V$ are endowed with the norm $\trinl \bullet \trinr_{\rm pw}\equiv\|\bullet\|_h|_{V+\M(\T)}$ and $a_h = a_{\rm pw} + {c}_{\rm P}$ for WOPSIP. %
\begin{lemma}[modified \eqref{eqn:H}] \label{lem:H_WOPSIP}
   Any  $(v_2,w_\M)\in P_2(\T)\times \M(\T)$ satisfies
\begin{align}\label{eqn:H_WOPSIP}
    a( J I_{\M} v_2, J w_{\M})- a_h(v_2,w_{\M})\le \| J \|
\trinl v_2- JI_\M v_2\trinr_{\pw} \trinl w_{\M} \trinr_{\pw}.
\end{align}
\end{lemma}
\begin{proof} For the WOPSIP scheme,  ${ c}_{\rm P}(v_2, w_\M) = 0$  for $(v_{2},v_{\M})\in P_2(\T)\times \M(\T)$ shows
	\begin{align*} %
a(JI_\M v_2, Jw_\M)	- a_{h}(v_2, w_\M)  &= 
a(JI_\M v_2, Jw_\M)-a_{\rm pw}(v_2, w_\M) \\
&=a_{\pw}(JI_\M v_2-v_2, Jw_\M) \le \|J\| 
\trinl v_2- JI_\M v_2\trinr_{\pw} \trinl w_{\M} \trinr_{\pw}
	\end{align*}
 with $a_{\rm pw}(v_2, w_\M- Jw_\M)=0 $ from \eqref{eqn:rightinverse} and the boundedness of $a_{\rm pw}(\bullet,\bullet)$ in the last line above.
\end{proof}
\noindent A careful revisit of the arguments in Section \ref{sec:Abstract a posteriori error control} with \eqref{eqn:H_WOPSIP} instead of \eqref{eqn:H} reveals a modified reliabiltiy
\begin{align}\label{eqn:reliability_WOPSIP}
    \| u- v_h\|_{h}  \lesssim \varrho  + \|(1-JI_\M)v_h\|_{h} +\|u_h-v_h\|_{h} %
\end{align}
with $\rho$ from \eqref{eqn:varrho} instead of $\rho(M)$ in Theorem \ref{thm:apost}. Theorem \ref{thm:apost_ns_final} can follow with Theorem \ref{thm:releff}.
\begin{proof}[Proof of Theorem \ref{thm:apost_ns_final}.b]
The a priori results for WOPSIP \cite[Thm.~8.13]{CCNNGRDS22} reveal the convergence under uniform mesh refinement in the norm $\|\bullet\|_{\rm P}$ and provide universal constants $\varepsilon_0,\delta_0>0$ such that, for any $\T\in\TT{\delta_0}$, a unique discrete solution $u_h\in V_h$ to \eqref{eqn:discrete} exists with $\|u-u_h\|_{\rm P}\leq \varepsilon_0$.
Since $\|\bullet\|_{\rm P}$ is stronger than $\|\bullet\|_h$ by \eqref{eqn:P_norm}, the convergence also follows in the weaker norm $\|\bullet\|_h$ and $h_{\rm max}\leq \delta_0$ implies the unique existence of a discrete solution $u_h\in V_h$ to \eqref{eqn:discrete} with $\|u-u_h\|_h\leq(1+\delta_0^2)^{1/2}\epsilon_0$. 
A density argument for $\varepsilon\coloneqq\min\{(1+\delta_0^2)^{1/2}\varepsilon_0,\beta/\left(3(1+\LP)\|\Gamma\|\right)\} $ leads to $\delta\leq\delta_0$ such that $\|u-u_h\|_h<\varepsilon$ for any $\T\in\TT\delta$.
This implies {\bf(L1)--(L3)} in the weaker norm $\|\bullet\|_h$ for $\delta,\varepsilon,\kappa\coloneqq2/3$, and $\varrho\coloneqq\beta/\left(3(1+\LP)\|\Gamma\|\right)$.
Hence, the setting of Section \ref{sec:Abstract a posteriori error control} applies and the proofs in Subsection~\ref{sec:apost} follow verbatim for $u_h,v_h \in V_h$, and $u \in V $ until \eqref{eqn:u-vh} that becomes
\begin{align} \label{eqn:triangle}
    \|u-v_h\|_{h} & \le 
     \|v_h -J I_\M v_h\|_{h}+
    \beta^{-1} (1-\kappa)^{-1} \|N(J I_\M v_h)\|_{V^*}
 \end{align}
 in the current setting.
The subsequent estimation of $\|N(J I_\M v_h)\|_{V^*}$ involves the split
$N(J I_\M v_h; y)=S_1+S_2+S_3$ as in \eqref{eqn:newb} for any $y \in Y \equiv H^2_0(\Omega)$ with $\trb{y}=1$ and the particular choice $y_h\coloneqq I_\M y\in \M(\T)$ in the definition of $S_1, S_3$ (instead of any $y_h\in P_2(\T)$ in Subsection~\ref{sub:varrho}).
The point of this modification is twofold. 
First, the supremum over $y\in Y$ with $w:=y- J I_\M y_h= y- J I_\M y \in V$ reveals
\begin{align}
    S_1:=a(J I_\M v_h,w) + \Gamma_{\pw}(Rv_h,Rv_h, w) -{F}(w) \le \rho.
\end{align}
Second, Lemma \ref{lem:H_WOPSIP} with $y_h\equiv I_\M y\in \M(\T)$ and $\trb{\bullet}_\pw\leq\|\bullet\|_h$ from \eqref{hnorm} result in \eqref{eqn:H_application} for $\Lambda_{\rm C}\coloneqq \|J\|$ without the need of \eqref{eqn:H}.
With this alternate derivation of \eqref{eqn:H_application}, the control of 
\begin{align}
    S_2&\coloneqq \Gamma_{\pw}(JI_\M v_h,JI_\M v_h,y) -\Gamma_{\pw}(JI_\M v_h,JI_\M v_h, y)\lesssim \|v_h -J I_\M v_h\|_{h},\\
    S_3&\coloneqq a(JI_\M v_h,JI_\M y)-F(JI_\M y) + \Gamma_{\pw}(Rv_h,Rv_h, JI_\M y)\lesssim \|v_h -J I_\M v_h\|_{h}+\|u_h-v_h\|_{h}\label{eqn:S3_bound_WOPSIP}
\end{align} %
follows verbatim from Subsection \ref{S2}--\ref{ssub:abstract_reliability_final} with $\|y_h\|_{h} = \trinl I_\M y\trinr_{\rm pw} \le \trinl y \trinr=1 \equiv M$.
This and \eqref{eqn:triangle}--\eqref{eqn:S3_bound_WOPSIP} verify the alternate reliability estimate \eqref{eqn:reliability_WOPSIP}. 
Theorem~\ref{thm:releff} further controls $\rho, \|(1-JI_\M)v_h\|_h$ in terms of the explicit a posteriori error terms $\mu(\T) +\eta(\T)\approx \sigma(\T)$  plus data oscillations ${\rm osc}_k(F,\T)$.
The efficiency of $\sigma(\T)$ with respect to $\|u-u_h\|_h$ plus data oscillations follows from Theorems~\ref{thm:releff}.b and~\ref{thm:apost}.b. This verifies \eqref{eqn:WOPSIP_NVS_a}. The sum of \eqref{eqn:WOPSIP_NVS_a} and $c_{\rm P}(v_h,v_h)=c_{\rm P}(u-v_h,u-v_h)$ with 
$ \|\bullet\|_h^2+c_{\rm P}=\trb{\bullet}^2+j_h+c_{\rm P}\approx \|\bullet\|_{\rm P}^2$
from \eqref{eqn:P_norm} reveal \eqref{eqn:WOPSIP_NVS_b} and conclude the proof.
\end{proof}
\section{Application to von K\'{a}rm\'{a}n equations }\label{sec:vke}
The von K\'{a}rm\'{a}n equations model a nonlinear plate \cite{CiarletPlates,Ciarlet2013} in two coupled PDE and require the product spaces
$X=Y=\bv := V \times V$ of $V\equiv H^2_0(\Omega)$
with norm  $\trinl\bullet\trinr$  
defined by $\trinl\boldsymbol\varphi\trinr:=(\trinl \varphi_1 \trinr^2+\trinl \varphi_2\trinr^2)^{1/2}$ for all $\boldsymbol\varphi=(\varphi_1,\varphi_2)\in \bv$. 
\subsection{Von K\'{a}rm\'{a}n equations}
\label{vke:model_intro}
Given any load $F\in H^{-2}(\Omega)$, the \vket~seek a solution ${\bf u} \equiv  ( \ua, \ub) \in  \bV $\;%
to
\begin{align}
\Delta^2 \ua =[\ua,\ub]+ F \;\;\text{ and }\;\;\Delta^2 \ub =-\half[\ua,\ua] 
\label{vkedGb}
\end{align}
in a bounded polygonal Lipschitz domain $\Omega \subset \mathbb{R}^2$.
Here and throughout this section, the (symmetric) von K\'{a}rm\'{a}n bracket $[\bullet,\bullet]$ reads $[\eta,\chi]:=\eta_{xx}\chi_{yy}+\eta_{yy}\chi_{xx}-2\eta_{xy}\chi_{xy}$.
 Recall the bilinear form ${a}_{\rm pw}$ and $\widehat V\equiv H^2(\T)$ from Subsection \ref{sub:triangulation, interpolation, and smoother}. %
 Let $\widehat{\bv}\coloneqq\widehat V\times \widehat V$ and define the trilinear forms $\Gamma_{\pw, 1}, \Gamma_{\pw, 2}:\widehat V\times \widehat V\times \widehat V$ and $\mathbf{\Gamma}_\pw:\widehat{\bv}\times \widehat{\bv}\times \widehat{\bv}$  by %
 \begin{align*}
     \Gamma_{\pw, 1}(\widehat{\xi},\widehat{\theta},\widehat{\varphi})&:=
-\sum_{T\in\T}\int_T [\widehat{\xi},\widehat{\theta}]\,\widehat{\varphi}\dx,\qquad%
\Gamma_{\pw, 2}(\widehat{\xi},\widehat{\theta},\widehat{\varphi}):=
\half\sum_{T\in\T}\int_T [\widehat{\xi},\widehat{\theta}]\,\widehat{\varphi}\dx,\\
\mathbf{\Gamma}_\pw(\widehat{\boldsymbol\xi },\widehat{\boldsymbol\theta},\widehat{\boldsymbol\varphi}) &:=\Gamma_{\pw,1}(\widehat{\xi_1},\widehat{\theta_2},\widehat{\varphi}_1)+\Gamma_{\pw,2}(\widehat{\xi_1},\widehat{\theta_1},\widehat{\varphi}_2)
 \end{align*}
for all $\widehat{\xi},\widehat{\theta},\widehat{\varphi} \in  \widehat V$ and $\widehat{\boldsymbol\xi }=(\widehat{\xi_1},\widehat{\xi_2}),\widehat{\boldsymbol\theta}=(\widehat{\theta_1},\widehat{\theta_2}),\widehat{\boldsymbol\varphi}=(\widehat{\varphi}_1,\widehat{\varphi}_2)\in \widehat{\bv}$.
Given $\displaystyle  a:= {a}_{\pw}|_{V \times V},$ %
$\Gamma_1:=\Gamma_{\pw,1}|_{V \times V\times V}$, and $\Gamma_2:=\Gamma_{\pw,2}|_{V \times V\times V}$, the weak formulation of~\eqref{vkedGb} seeks  ${\bf u} \equiv  ( \ua, \ub) \in  \bV $ with
\begin{align} \label{vk_weak}
a(\ua,\varphi)+ \Gamma_1(\ua,\ub,\varphi) = F(\varphi) \quad\text{and}\quad
a(\ub,\varphi) + \Gamma_2(\ua,\ua,\varphi)   = 0\quad\text{for all }\varphi\in V.       
\end{align} 

\medskip \noindent   For all $\boldsymbol\theta=(\theta_1,\theta_2)$ and $\boldsymbol\varphi=(\varphi_1,\varphi_2)\in  \bv$, define %
$${\boldsymbol a}(\boldsymbol\theta,\boldsymbol\varphi) :={} a(\theta_1,\varphi_1) + a(\theta_2,\varphi_2),\quad\text{and}\quad\mathbf F(\boldsymbol\varphi) \coloneqq F(\varphi_{1})
.$$ 
Given $\mathbf{\Gamma}\coloneqq \mathbf{\Gamma}_\pw|_{V\times V\times V}$, the vectorised formulation of~\eqref{vk_weak} seeks ${\bf u}=(\ua,\ub)\in \bv$ such that
\begin{equation}\label{VKE_weak}
\mathbf N({\bf u};\boldsymbol\varphi):={\boldsymbol a}({\bf u},\boldsymbol\varphi)+\mathbf \Gamma({\bf u},{\bf u},\boldsymbol\varphi)- \mathbf{F}(\boldsymbol\varphi)=0\fl \boldsymbol\varphi \in \bv.
\end{equation}
Given any $F\in V^*\equiv H^{-2}(\Omega)$, there exist solutions to \eqref{vk_weak}, which are possibly singular but carry elliptic regularity; the weak solution is unique and a regular root in the case of small loads, cf.~\cite{MN:ConformingFiniteElement2016,Ciarlet2013,Kni:ExistenceTheoremKarman1967} for proofs. The a posteriori error analysis below concerns some approximation ${\bf v}_{h}\in \bV_h \coloneqq  V_h\times V_h$  to a regular root  ${\bf u} \in \bV$ of the continuous problem~\eqref{vk_weak}.
}
\subsection{A posteriori error control for five quadratic discretizations} \label{sec:vke_disc}
This subsection applies the abstract a posteriori error analysis from Section \ref{subsection:new}--\ref{sec:nse} to the Morley, dGI\&II, $C^0$IP, and WOPSIP schemes for~\eqref{VKE_weak}.
Recall the discrete space $ V_h$ from Subsection~\ref{sec:discrete_schemes} together with the bilinear forms $\mathsf{b}_h$ and $\mathsf{c}_h$
from Table \ref{tab:spaces} for the five methods.
For any $\boldsymbol\theta_h\equiv(\theta_{h,1},\theta_{h,2}),\boldsymbol\varphi_h\equiv(\varphi_{h,1},\varphi_{h,2})\in X_h\equiv Y_h\equiv\bv_h\equiv V_h\times V_h$, the discrete bilinear form ${\boldsymbol a}_h : \bv_h \times \bv_h \rightarrow \mathbb{R}$ reads
\begin{equation}\label{eqn:apw_vke}
    \begin{aligned}
         {\boldsymbol a}_{h}(\boldsymbol\theta_h,\boldsymbol\varphi_h)
:={}& a_{\rm pw}(\theta_{h,1},\varphi_{h,1})+ \mathsf{b}_h(\theta_{h,1},\varphi_{h,1}) +  \mathsf{c}_h(\theta_{h,1},\varphi_{h,1})  \\
&+ a_{\rm pw}(\theta_{h,2},\varphi_{h,2}) +  \mathsf{b}_h(\theta_{h,2},\varphi_{h,2}) + \mathsf{c}_h(\theta_{h,2},\varphi_{h,2}).
    \end{aligned}
\end{equation}
The second dG scheme replaces $a_{\rm pw}$ by $(\Delta_{\pw} \bullet, \Delta_{\pw} \bullet)_{L^2(\Omega)}$.
Let  $\boldsymbol{R} \in \{\mathbf{id}, \boldsymbol{I}_{\rm M}, \boldsymbol{JI}_{\rm M} \}$ and $\boldsymbol{P}\equiv\boldsymbol Q\equiv\boldsymbol S \equiv\boldsymbol J\boldsymbol I_\M$ denote the vectorized versions of the respective operators from Subsection \ref{sub:triangulation, interpolation, and smoother} that apply componentwise.
The discrete scheme for~\eqref{VKE_weak} 
 seeks a solution ${\bf u}_{h}\in \bV_h$ to
\begin{align}\label{eqn:gen_scheme_vke}
\hspace{-0.7cm}\boldsymbol{N}_h({\bf u}_h;\Phi_h) := \boldsymbol a_{h}({\bf u}_{h},\boldsymbol\varphi_h)
+\mathbf\Gamma_{\pw}(\boldsymbol {R}{\bf u}_{h},\boldsymbol{R}{\bf u}_{h},
\boldsymbol{S}\boldsymbol\varphi_h)-\mathbf F( \boldsymbol{ Q}{\boldsymbol\varphi}_{h})=0 \fl \boldsymbol\varphi_h \in  \bv_h.
\end{align}
%
%
%
%
%
%
%
%
%
%
%
%
%
%
%
%
%
%
\begin{comment}
\begin{table}[H]
	\centering
	\begin{tabular}{|c|c|c|c|c|c|}
		\hline
		Scheme     & Morley & dG I       & $C^0$IP     & WOPSIP & dG II \\ \hline
		%
		$X_h = Y_h = \bV_h$          &  $\mathbf{M}(\cT)$      & $\boldsymbol{P}_2(\cT)$      & $\boldsymbol{S}^2_0(\cT)$ & $\boldsymbol{P}_2(\cT)$  & $\boldsymbol{P}_2(\cT)$      \\ \hline
		\begin{minipage}{2cm}
			\centering
			$\widehat{X} = \widehat{Y} = \widehat{\bV}  =\bV + \bV_h$
		\end{minipage}   &   $\bV + \mathbf{M}(\cT)$     & $ \bV + \Stb$          & $ \bV + \boldsymbol{S}^2_0(\cT)$  & $ \bV +  \Stb$   &  $ \bV + \Stb$     \\ \hline
		$\|\bullet \|_{\widehat{X}}$            & $\trinl \bullet \trinr_{\rm pw}$     & $\| \bullet \|_{\rm dG}$       &    $\| \bullet \|_{\rm IP}$  &  $\| \bullet \|_{\rm P}$ &$\| \bullet \|_{\rm dG}$      \\ \hline
%
%
		$I_{h}$ & ${\mathbf {id}}$ & ${\mathbf {id}}$ &$\boldsymbol{I}_{\rm C}$  & ${\mathbf  {id}}$ & ${\mathbf {id}}$ \\ \hline 
			$I_{X_h}=I_{\bV_h} = I_{h}\boldsymbol{I}_{\rm M}$           &  $\boldsymbol{I}_{\rm M}$ & $\boldsymbol{I}_{\rm M}$ & $\boldsymbol{I}_{\rm C}\boldsymbol{I}_{\rm M}$ & $\boldsymbol{I}_{\rm M}$ &  $\boldsymbol{I}_{\rm M}$    \\ \hline
%
%
%
%
%
%
	\end{tabular}
	\caption{Spaces, operators, and norms in Section \ref{sec:vke}.}
	\label{tab:spaces_vke}
\end{table}
\end{comment}
%
Recall $\tau_{\E},\nu_{\E},[\bullet]_{\E}$, and the piecewise integral mean operator $\Pi_{\E,0}$ from Subsection~\ref{sub:triangulation, interpolation, and smoother}. 
Set $\vartheta=1$ for Morley, dG I\& II, WOPSIP, and
$\vartheta=0$ for C$^0$IP.
The a posteriori error estimator $\boldsymbol{\sigma}(\T) := \sqrt{\sum_{T \in \T} \boldsymbol{\sigma}^2(T)}$ for some approximation ${\bf v}_h= (\vha,\vhb)\in \bV_h$  to a regular root  ${\bf u} \in \bV$ of the continuous problem \eqref{vkedGb} has on $T\in\T$ the contribution
\begin{equation*}%
    \begin{aligned}
\boldsymbol{\sigma}^2(T)& := |T|^2 \Big( \|\Pi_k f_0+[\vha,\vhb] - {\rm div}_{\pw} (\Pi_k f_1) +{\rm div}^2_{\pw}  (\Pi_k f_2) \|^2_{L^2(T)} + \|[\vha,\vha] \|^2_{L^2(T)} \Big) \\
& \quad + |T|^{3/2} \|\jump{\Pi_k f_1 - {\rm div} (\Pi_k f_2) -
\partial_s ( \Pi_k f_2) \tau_{\E}}  \cdot \nu_{\E} \|^2_{L^2(\partial T \setminus \partial \Omega)}\\
& \quad + |T|^{1/2} \Big( \big\|[D^2\vha]_{\E}\tau_{\E}\big\|_{L^2(\partial T)}^2+ { \big\| (1-\vartheta\Pi_{\E,0}) \big[\big(\Pi_k f_2-  D^2 \vha\big)\nu_{\E}\big]_{\E} \cdot  \nu_{\E} \big\|_{L^2(\partial T \setminus \partial \Omega)}^2}
\\
& \quad \phantom{|T|^{1/2} \Big(}  +   \big\|[D^2\vhb]_{\E}\tau_{\E}\big\|_{L^2(\partial T)}^2+ { \big\| (1-\vartheta\Pi_{\E,0}) \big[{\big( D^2 \vhb\big)\nu_{\E}}\big]_{\E} \cdot  \nu_{\E} \big\|_{L^2(\partial T \setminus \partial \Omega)}^2} \Big)
\\
& \quad + \sum_{\ell=1,2} \sum_{E \in \mathcal{E}(T)} \Big(\big| \Pi_{E,0} \big[\partial_\nu \vhl\big]_E \big|^2 + |T|^{-1}\hspace{-0.6em} \sum_{z \in \mathcal{V}(E)} \big|\big[\vhl]\big]_E(z)\big|^2 \Big).  
    \end{aligned}
\end{equation*}
Abbreviate $\mathbf{c}_{\rm P}({\bf v}_h,{\bf v}_h)\coloneqq c_{\rm P}(\vha,\vha)+c_{\rm P}(\vhb,\vhb)$ with the WOPSIP penalty form $c_{\rm P}$ from \eqref{eqn:c_P_def}.%
\begin{theorem}[a posteriori error control]%
\label{thm:apost_vke_final}
Given a regular root ${\bf u} \in \bV$  to \eqref{VKE_weak} with $F\in H^{-2}(\Omega)$ and $k\in\mathbb N_0$, there exist $\varepsilon, \delta,\varrho 
>0$ such that \eqref{eqn:gen_scheme_vke} has a unique discrete solution ${\bf u}_h\in \bV_h$ to \eqref{eqn:gen_scheme_vke} with $\|  {\bf u}-{\bf u}_h\|_{h} \le\varepsilon$ for any $\displaystyle\cT\in\bT(\delta)$ and the following holds for $F$ given as in \eqref{eqn:newone}.
$(a)$ For the Morley, dG I \& II, and $C^0$IP scheme and $\boldsymbol{P}\equiv \boldsymbol {Q} \equiv \boldsymbol{S}\equiv \boldsymbol {J} \boldsymbol{ I}_\M$, any ${\bf v}_h \in \bV_h$ with $\|{\bf u}_h-{\bf v}_h\|_h \le \varrho$ satisfies
$$\| {\bf u}-{\bf v}_h \|_{h}  + {\rm osc}_k(F,\T) \approx\boldsymbol{\sigma}(\T)  + \|{\bf u}_h-{\bf v}_h\|_{h} + {\rm osc}_k(F,\T).$$
$(b)$ For the WOPSIP scheme and $\boldsymbol{P}\equiv \boldsymbol {Q} \equiv \boldsymbol{S}\equiv \boldsymbol {J} \boldsymbol{ I}_\M$, any ${\bf v}_h \in \bV_h$ with $\|{\bf u}_h-{\bf v}_h\|_h \le \varrho$ satifsies
\begin{align*}%
\| {\bf u}-{\bf v}_h \|_{h}  + \|{\bf u}_h-{\bf v}_h\|_h + {\rm osc}_k(F,\T) &\approx\boldsymbol{\sigma}(\T)  + \|{\bf u}_h-{\bf v}_h\|_{h} + {\rm osc}_k(F,\T),\\
\| {\bf u}-{\bf v}_h \|_{\rm P}  + \|{\bf u}_h-{\bf v}_h\|_h + {\rm osc}_k(F,\T) &\approx\boldsymbol{\sigma}(\T) + \mathbf{c}_{\rm P}({\bf v}_h,{\bf v}_h)^{1/2} + \|{\bf u}_h-{\bf v}_h\|_{h} + {\rm osc}_k(F,\T).
\end{align*}
\end{theorem}
\begin{proof}
The proof employs the a priori analysis in \cite[Sec.~9]{CCNNGRDS22} for the existence of a local unique discrete solution and follows the lines of that of Theorem \ref{thm:apost_ns_final} for all the components behind the vector notation of this section. Indeed,
Theorem \ref{thm:apost}.a and Theorem \ref{thm:releff}.a provide $\| {\bf u}-{\bf v}_h \|_{h}\lesssim \boldsymbol{\varrho}+\|{\bf v}_h-\boldsymbol{JI}_\M{\bf v}_h\|_h+ \|{\bf u}_h-{\bf v}_h\|_{h}$, where $\boldsymbol{\rho}\equiv \rho_1+\rho_2$ from \eqref{eqn:varrho} splits into the components
\begin{align}
        \rho_1&\coloneqq
\sup_{\stackrel{y \in Y, \|y \|_Y = 1}{z:=y- J I_\M I_hI_\M y} } 
\big( F(z)-a( J I_\M \vha,z) -
\Gamma_{\pw,1}(R\vha, R\vhb, z)  \big),\label{eqn:varrho_1}\\
        \rho_2&\coloneqq
\sup_{\stackrel{y \in Y, \|y \|_Y = 1}{z:=y- J I_\M I_hI_\M y} } 
\big(-a( J I_\M \vhb,z) -
\Gamma_{\pw,2}(R\vha, R\vha, z)  \big).\label{eqn:varrho_2}
\end{align}
The control of $\rho_1,\rho_2$ as in Theorem \ref{thm:releff} amounts in a large number of terms gathered together in the estimator $\boldsymbol{\sigma}(\T)$. As there is no additional mathematical difficulty, the further details are omitted.
\end{proof}
\subsection{Single force} 
\label{sub:Single_forces}
Practical plate problems also concern line loads and single forces as discussed in \cite{CCBGNN22} for the linear biharmonic plate.
There are two amazing observations regarding a single force $\lambda_{\zeta} \pl$ with strength $\lambda_{\zeta} \in {\mathbb R}$ and the Dirac delta distribution $\pl \in H^{-2}(\Omega)$ at a finite set $A\subset \Omega$ of points $\zeta \in A$. 
First, single loads $\lambda_{\zeta} \pl$ decouple in the a posteriori error analysis %
from a general source $F \in H^{-2}(\Omega)$ given~by
\begin{align}\label{eqn:F_extended}
    F (\varphi) := \int_{\Omega} (f_0\; \varphi + f_1 \cdot \nabla \varphi + f_2 : D^2 \varphi) \dx + \sum_{\zeta\in A}\lambda_{\zeta} \pl(\varphi)\quad \text{for all } \varphi \in H^2_0(\Omega).
\end{align}
 in terms of 
 Lebesgue functions $f_0 \in L^2(\Omega), f_1 \in L^2(\Omega; {\mathbb R}^2), f_2 \in L^2(\Omega; {\mathbb S})$ as in \eqref{eqn:newone}.
 Second, a point load at an interior vertex $\zeta \in {\mathcal V}(\Omega)$ leads to a load in the discrete problem \eqref{eqn:discrete} but has no contribution to the a posteriori error estimate because the test functions for $F$ in $\rho_1$ from \eqref{eqn:varrho_1} vanish at all vertices $\mathcal V$.
 This has already been observed in \cite{ccnn2021, CCBGNN22} for linear problems. 
 Hence we consider a finite family $(\lambda_\zeta\pl: \zeta\in A)$ of single forces at $A\subset\Omega, |A|<\infty$ and distinguish $A\cap\mathcal V$ with no contributions $\mu(\zeta)=0$ and $A\setminus\mathcal V$ with a contribution $\mu(\zeta)$ to the a posteriori error control as follows.
 
\begin{comment}
 %
This subsection focuses on one single force $F=\lambda \pl$ with strength $\lambda \in {\mathbb R}$ and the Dirac delta distribution $\pl \in H^{-2}(\Omega)$ at a point $\zeta \in \Omega$. There are two amazing observations. First, the term $\lambda \pl$ (even a finite sum of those) decouples from $F= \lambda \pl+G\in H^{-2}(\Omega)$ with $G$ given in \eqref{eqn:newone} into  Lebesgue functions $f_0,f_1,f_2$ and hence all those functions are set zero in this subsection.
%
 Second, for a vertex $\zeta \in {\mathcal V}(\Omega)$, the functional $\pl$ applies to a test function $z= w - JI_\M I_h I_\M w$ that vanishes at ${\mathcal V} $ and so $\lambda(\pl )(z)= \lambda \: z(\zeta) =0$. Hence there remains no extra contribution in the a posteriori error estimate. In other words, $\zeta \in {\mathcal V}(\Omega)$ leads to a load in the discrete problem \eqref{eqn:discrete} but is not visible in the error estimate; this is already observed for linear problems in \cite{ccnn2021, CCBGNN22}. 
    
\end{comment}
\newcommand{\etaz}{\mu(\zeta)}
 \medskip \noindent
 Consider a single force $\lambda_\zeta\pl $ at a generic position $\zeta \in \Omega \setminus {\mathcal V}$ that is different from a vertex of the triangulation so that at most two triangles $T \in {\mathcal T}(\zeta):= \{ K \in \T: \zeta \in K\}$ contain $\zeta \in T$.
 Let  $h_\zeta\coloneqq\min {\{ }|T|^{1/2}: T \in \T(\zeta) {\}}, \omega(\zeta)\coloneqq\mathrm{int}(\cup \T(\zeta))$ and suppose the separation assumption
 \begin{align}\label{eqn:single_load_assumption}
     {\rm dist}(\zeta, {\mathcal V}) \gtrsim h_\zeta\quad\text{and}\quad |A\cap\omega(\zeta)|=1
     \qquad\text{for all }\zeta\in A\setminus\mathcal V.
 \end{align}
 The following result extends Theorem \ref{thm:apost_vke_final} to right-hand sides with single forces at $A$.
 Define the novel estimator $\etaz\coloneqq|\lambda_\zeta| h_\zeta$ for a single load at $\zeta\in A$ and set $\mu(\mathcal M)\coloneqq\sum_{\zeta\in \mathcal M}\etaz$ for all $\mathcal M\subset A$.
 \begin{theorem}[single forces]\label{thm:apost_vke_dirac}
     Under the assumptions of Theorem \ref{thm:apost_vke_final}.a for $F\in H^{-2}(\Omega)$ given in \eqref{eqn:F_extended} for single forces at a finite set $A\subset\Omega$ with \eqref{eqn:single_load_assumption}, the Morley, dG I \& II, and $C^0$IP scheme satisfy
     $$\| {\bf u}-{\bf v}_h \|_{h}  + {\rm osc}_k(F,\T) \approx\boldsymbol{\sigma}(\T) + \mu(A\setminus\mathcal V) + \|{\bf u}_h-{\bf v}_h\|_{h} + {\rm osc}_k(F,\T).$$
     An analog to Theorem~\ref{thm:apost_vke_final}.b holds for the WOPSIP scheme where $\mu(A\setminus\mathcal V)$ is added on the respective right-hand sides.
 \end{theorem}

\begin{example}\label{ex:separation_condition}
    This separation condition \eqref{eqn:single_load_assumption} is met in the numerical benchmark in Subsection~\ref{sub:Numerical results for the Von Karman equation} for the centroid $\zeta=(-1/6;-1/6)$ of the L-shaped domain $\Omega$ that lies on an edge $E=\mathrm{conv}\{A,B\}\in\E$ parallel to the main diagonal for all triangulations (occuring from newest-vertex bisections) with $\zeta=1/3A+2/3B$ and ${\rm dist}(\zeta, {\mathcal V}) = h_E/3\geq \sqrt2h_\zeta/3 $ for the right-isosceles triangles in Subsection~\ref{sub:Numerical results for the Von Karman equation}.
\end{example}
\begin{proof}[Proof of Theorem \ref{thm:apost_vke_dirac}] The modifications to the source $F$ only enter the estimation of the term $\rho_1$ from \eqref{eqn:varrho_1} in Theorem \ref{thm:apost_vke_final}.
The reliability follows 
from Section~\ref{subsection:new} plus the analysis, for $\zeta\in A$, of the extra terms 
\[ |\lambda_{\zeta} \delta_{\zeta} (z)| \equiv |\lambda_{\zeta} z(\zeta)| \le \begin{cases}
    |\lambda_{\zeta}|\,C_{\rm BH}h_\zeta   |z|_{H^2(T)}&\text{if }\in A\setminus\mathcal V,\\
    0&\text{if }\in A\cap \mathcal V
\end{cases} \]
for any $z\coloneqq w-JI_\M I_hI_\M w$ with $w\in V$.
Indeed, $z$ vanishes at ${\mathcal V}$ and a Bramble-Hilbert lemma scales the Sobolev embedding $H^2(T_{\rm ref}) \hookrightarrow C(T_{\rm ref})$ from a reference triangle $T_{\rm ref}$ to $T$ with a constant $C_{\rm BH}$ that exclusively depends on the shape-regularity of $T$. Recall \eqref{eqn:C_stability} for $\trinl z \trinr \le C_{\rm stab} \trinl w \trinr$ and $\trinl w \trinr =1$ to deduce 
$\lambda_{\zeta} \delta_{\zeta} (z) \le C_{\rm BH} C_{\rm stab} \,\mu(\{\zeta\}\setminus \mathcal V)$
and involve this estimate in Subsection~\ref{appl:rel}. This outlines the proof of $\rho_1 \lesssim \boldsymbol{\sigma}(\T) + \mu(A\setminus\mathcal V) + {\rm osc}_k(F,\T)$;
the remaining details are straightforward from Theorem \ref{thm:apost_vke_final} and hence omitted. The efficiency of the additional a posteriori error terms $\mu(\zeta)$ for any $\zeta\in A\setminus\mathcal V$ requires the design of a test function $\psi \in V$ with a list of properties: 
\begin{equation}\label{eqn:psi_def}
    \begin{aligned}
        & \psi(\zeta)=1, \; \psi=0 \text{ at } {\mathcal V}, \;
{\rm supp}\: \psi \subseteq \overline{\omega(\zeta)}, \;
\psi|_T \perp P_k(T) \text{ in } L^2(T) \text{ for any triangle }
T \in \T,
\\
&  \; \psi|_E \perp P_k(E) \text{ in } L^2(E), 
\text{ and }
\nabla \psi|_E  \perp (P_k(E))^2 \text{ in } L^2(E)^2 \text{ along any edge } E \in \E.
    \end{aligned}
\end{equation}
This function can always be constructed and Supplement B provides an elementary design of $\psi$ in terms of Jacobi polynomials if $\zeta\in E\in\E$. An important detail is the scaling $\trinl \psi \trinr \approx h_{\zeta}^{-1}$ that requires the separation condition ${\rm dist}(\zeta, {\mathcal V}) \gtrsim h_\zeta $, while
 the universal case with $\zeta\in\Omega$ arbitrary involves a more refined analysis with a weight that is left for future research.
Recall $\Pi_kF$ from \eqref{eqn:Pi_F_def} and consider $\Lambda\coloneqq \Pi_kF-a_\pw(\vha,\bullet)-\Gamma_{\pw, 1}(\vha,\vhb, \bullet)$ as in
\eqref{eqn:newthree} with 
 \begin{align*}
    {\Lambda}_0:= \Pi_k f_0 + [\vha,\vhb] \in P_k(\T),\;\;
\Lambda_1:= \Pi_k f_1 \in  P_k(\cT; {\mathbb R}^2), \;\;
\Lambda_2 := \Pi_k f_2 - D^2_{\pw} \vha \in P_k(\cT; {\mathbb S})
\end{align*}
 The many orthogonalities \eqref{eqn:psi_def} of $\psi$ enter the final stage in two piecewise integration by parts for 
\begin{align}\nonumber
 \Lambda(\psi)
  \equiv  \int_{\Omega} \big( \Lambda_0 \psi + \Lambda_1 \!\cdot \!\nabla \psi + \Lambda_2\!:\! D^2 \psi \big) \dx &=
  \sum_{E \in \E} \int_E \big( \psi \jump{\Lambda_1-\div_{\pw} \Lambda_2}_E\! \cdot\! \nu_E  + \nabla \psi\! \cdot \!\jump{\Lambda_2}_E \nu_E \big)\,\ds\\
  \label{eqn:algebraicidentity}
 &\; + \sum_{T \in \T} \int_T \big(  \Lambda_0 - \div_\pw \Lambda_1 + 
\div_\pw^2 \Lambda_2 \big)\, \psi\dx =0.
\end{align}
In fact, $\psi|_T \perp \Lambda_0- \div_{\pw} \Lambda_1+ \div^2_{\pw} \Lambda_2 \in P_k(T)$ in $L^2(T)$ for $T \in \T$, $\psi|_E \perp \jump{\Lambda_1-\div_{\pw} \Lambda_2}_E \in P_k(E) $ in $L^2(E)$ as well as $(\nabla \psi)|_E \perp \jump{\Lambda_2}_E \nu_E  \in P_k(E)^2 $ in $L^2(E)^2$ for all $E \in \E$ make each of the integrals vanish individually. 
Since $\zeta$ is the only element of $A\cap\omega(\zeta)$ by \eqref{eqn:single_load_assumption}, the properties of $\psi$ imply $\sum_{a\in A}\lambda_{a}\delta_{a}(\psi) = \lambda_{\zeta}$.
This, the algebraic identity \eqref{eqn:algebraicidentity}, and the first component $a(\ua,\psi) + \Gamma_1(\ua,\ub, \psi) = F(\psi)$ of the problem \eqref{vk_weak} (and~\eqref{VKE_weak}) with exact solution ${\bf u}=(\ua,\ub) \in \bV$ lead to 
\[ \lambda_{\zeta}= %
a_\pw(\ua-\vha,\psi)
+ \Gamma_{\pw,1}(\ua-\vha,\ub,\psi) + \Gamma_{\pw,1}(\vha, \ub-\vhb,\psi) + \Big(\Pi_kF+\sum_{a\in A}\lambda_{a}\delta_{a}-F\Big)(\psi). \]
A Cauchy inequality for $a_\pw$, the boundedness of $\Gamma_{\pw,1}$, and a routine estimation of the last term  as in \eqref{eqn:fbound} with the scaling $\|\psi\|_{H^s(T)}\leq h_{T}^{2-s}C_{\rm BH}\|\psi\|_{H^2(T)}$ for $s=0,1,2$ and $T\in\T(\zeta)$ from a Bramble-Hilbert lemma provide
\[ h_\zeta^{-1}  \etaz =|\lambda_\zeta| \le \trinl \psi \trinr
\Big(\big(1+\|\Gamma_{\pw,1}\|( \trinl \ub\trinr + \| \vha \|_h )\big) \| {\bf u}-{\bf v}_h \|_h + C_{\rm BH}{\rm osc}_k(F,\T)\Big). \]
The scaling $\trinl \psi \trinr \approx h_\zeta^{-1} $ and the boundedness of $\trinl \ub\trinr + \| \vha \|_h\leq \trb{{\bf u}}+\|{\bf v}_h\|_h$ by \eqref{eqn:vhbound} imply $\etaz \lesssim \| {\bf u}-{\bf v}_h \|_h+{\rm osc}_k(F,\T)$.
Since the set $|A\setminus\mathcal V|\lesssim 1$ is finite, the sum over all $\zeta\in A\setminus\mathcal V$ concludes the efficiency of $\mu(A\setminus\mathcal V) \lesssim \| {\bf u}-{\bf v}_h \|_h+{\rm osc}_k(F,\T)$.
Theorem \ref{thm:releff}.c and the quasi-optimality \eqref{quasioptimalsmootherP} of the smoother $P\equiv JI_\M$ provide the efficiency of  $\|[D^2v_h]_{\E}\tau_{\E}\|_{L^2(\E)}+j_h(v_h,v_h)^{1/2}$ for $v_h=\vha,\vhb$.
The efficiency of the other terms in $\boldsymbol{\sigma}(\T)$ does not follow verbatim, but a correction of standard (cubic volume and quadratic edge) bubble functions by a multiple of $\psi$ so that the resulting sum vanishes at $A$ decouples the contributions and leads to local efficiency as in \cite[Sec.~7.4]{CCBGNN22}.
\end{proof}

\section{Numerical experiments}%
\begin{figure}[]
	\centering
	\hbox{
	\hspace{1.5em}\includegraphics{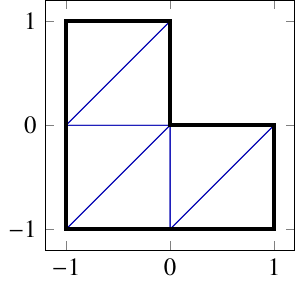}
	\hspace{1.5em}\includegraphics{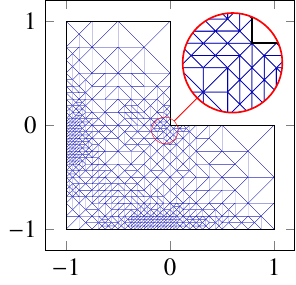}
	\hspace{1.5em}\includegraphics{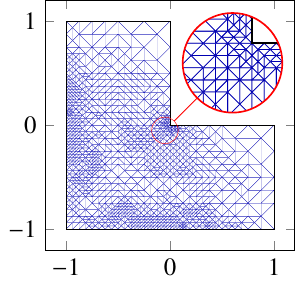}
	}
	\caption{Initial (left) and adaptive triangulations with $|\T|=1083$ (middle) and $|\T|=2044$ (right) triangles for
		the
	Morley FEM in Subsection \ref{sub:Singular solution on the L-shaped domain}}
	\label{fig:NVS_Grisvard_mesh}
\end{figure}

\label{sec:Numerical experiments}
This section compares the uniform and adaptive Morley, $C^0$IP, and dG FEM for the Navier-Stokes and von
K\'arm\'an equations in 2D on triangulations of the L-shaped domain of Figure \ref{fig:NVS_Grisvard_mesh} and~\ref{fig:VKE_Dirac_meshes} into triangles.
\subsection{Numerical realization}%
	\label{sub:Numerical realization}
 The Newton scheme allows the approximation $v_h$ up to machine precision of a root $v_h=u_h$ to \eqref{eqn:discrete} and so we disregard the algebraic error $\|u_h-v_h\|_h=0$. Supplement~C provides algorithmic details on the implementation of the nested iterations and the termination criterion. 
 This section presents numerical evidence on the theoretical results for the a posteriori error estimators $\sigma(\T)$ from Sections~\ref{sec:nse}--\ref{sec:vke} and the related standard D\"{o}rfler marking adaptive algorithm with newest-vertex bisection. 

\subsection{Navier-Stokes equations on L-shaped domain}%
\begin{figure}[]
		\centering
		\hspace{-3em}\includegraphics{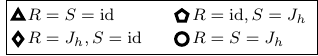}\hspace{7em}
		\includegraphics{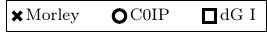}\\
		\hspace*{-2em}\hbox{
			\includegraphics{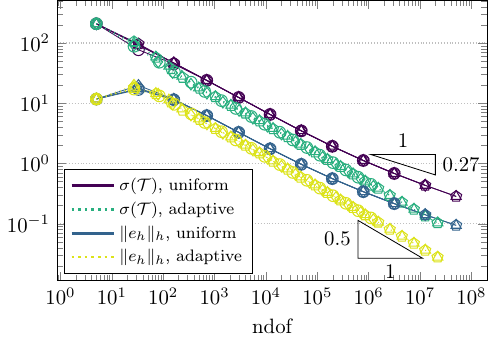}
			\includegraphics{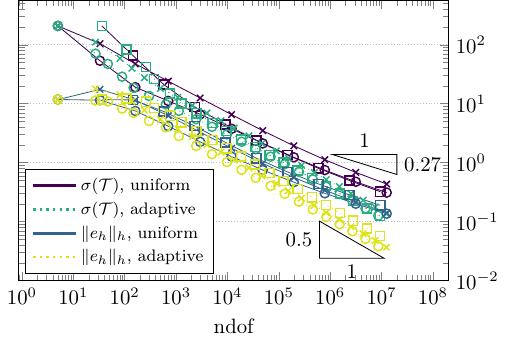}
		}
		\caption{Convergence history plot of the error $\|e_h\|_h$ and the estimator $\sigma(\T)$ for Morley FEM with different choices of $R, S\in\{{\rm id}, J I_\M\}$ (left) and Morley,
			C0IP, dG I FEM (right) for the singular solution $u$ in Subsection \ref{sub:Singular solution on the L-shaped domain}
		 }
		\label{fig:Square_exact_h}
	\end{figure}

\label{sub:Singular solution on the L-shaped domain}
 The \textit{singular solution} from Grisvard \cite{grisvard_singularities_1992} for the L-shaped domain $\Omega=(-1,1)^2\setminus[0,1)^2$ reads
	\begin{align*}
		u(r, \varphi)
			&=\big(r^2\sin(\varphi)^2-1\big)^2\big(r^2\cos(\varphi)^2-1\big)^2r^{1+\mu}\mu^2\;\xi\big(\varphi)
	\end{align*}
   in polar coordinates with interior angle $\omega=3\pi/2$ at the origin, $\mu=0.54448$, and a smooth function $\xi$ given \cite[Eqn.~3.2.9]{grisvard_singularities_1992} (therein denoted as $\xi(\varphi) = u(\mu,\varphi-\pi/2)$).
	This function $u \in H^2_0(\Omega) \cap
 H^{2+\sigma}(\Omega)$ for $\sigma < \mu $ serves as an exact solution to 
 \eqref{NS_eqa} with computed source term $F\equiv f \in L^2(\Omega)$.
 Figure \ref{fig:Square_exact_h} displays the expected suboptimal experimental convergence rate $\sigma/2$ on uniformly refined triangulations.
 The a posteriori error analysis in Section \ref{sec:nse} motivates the standard adaptive algorithm driven by the refinement indicators
 $\sigma^2(T)$ for a triangle $T\in\T$ equal to
 \begin{align} \label{eqn:eta1}
		&|T|^2\|f\|_{L^2(T)}^2  +|T|^{3/2}\|[\Delta
		u_h\mathrm{Curl} u_h]_{\E}\cdot\nu_{\E}\|_{L^2(\partial T \setminus \partial \Omega)}^2+\vartheta |T|^{1/2}
  \|[\partial_{\nu\nu}u_h]_{\E}\|_{L^2(\partial T \setminus \partial \Omega)}^2\nonumber\\
  +&|T|^{1/2} \big\|[D^2u_h]_E\tau_E\big\|_{L^2(\partial T)}^2  +\sum_{E \in \mathcal{E}(T)} \Big(\big| \Pi_{E,0} \jump{\partial_\nu u_h}_E \big|^2 + |T|^{-1}\hspace{-0.6em} \sum_{z \in \mathcal{V}(E)} \big|\jump{u_h}_E(z)\big|^2 \Big)
	\end{align}
 with $\vartheta=1$ for $C^0$IP and $\vartheta=0$ otherwise. 
 For all choices of the operators $R,S\in\{{\rm id}, I_\M, J I_\M\}$ shown for Morley on the left in Figure~\ref{fig:Square_exact_h}, the adaptive algorithm recovers optimal convergence rates of the error $e_h\coloneqq u-u_h$ in the norm $\|\bullet\|_h$ and, as implied by Theorem~\ref{thm:apost_ns_final}, the error estimator $\sigma(\T)\coloneqq\sqrt{\sum_{T\in\T}\sigma^2(T)}$.
 The first competition of the lowest-order Morley, dGI, and $C^0$IP scheme with parameters $\sigma_{\ip}\coloneqq \sigma_1\coloneqq \sigma_2\coloneqq 20$ in \eqref{eqn:cdG_def} and $\theta=1$ in Table \ref{tab:spaces} reveals an overall comparable performance with the smallest error for given number of dofs from the $C^0$IP scheme shown in Figure~\ref{fig:Square_exact_h}.
 The undisplayed efficiency indices $EF\coloneqq \sigma(\T)/\|e_h\|_h$ range between $1.5$ and $4$ on meshes with at least $1000$ dof.
 Figure \ref{fig:NVS_Grisvard_mesh} displays the initial triangulation and a typical output of the adaptive algorithm with the expected local refinement towards
	the singularity at the origin.
	The additional mild refinement near the sides opposite to the origin appears for all schemes with different intensity and is interpreted as
	a boundary layer already observed for the linear biharmonic problem, e.g., in \cite{CCDGJH14}.

\subsection{Von K\'arm\'an problem with a point load}%
\begin{figure}[]
				\centering
		\hbox{
					\hspace{0em}\includegraphics{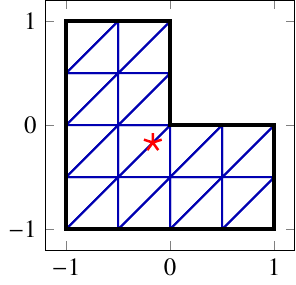}
					\hspace{1em}\includegraphics{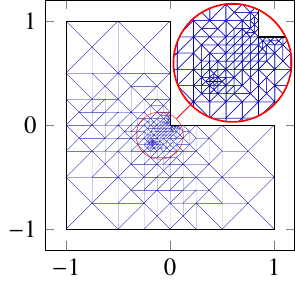}
				\hspace{1em}\includegraphics{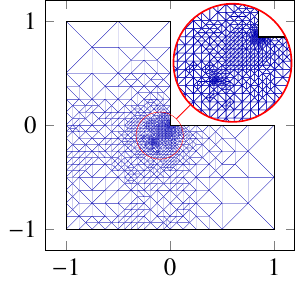}
				}
				\caption{Initial mesh for the von K\'arm\'an problem
    with the centroid \includegraphics{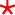} (left) and adaptive triangulations with $|\T|=416$ (middle) and $|\T|=1678$ (right) triangles from
			the
		Morley FEM for $R=S=\id$}
				\label{fig:VKE_Dirac_meshes}
		\end{figure}

\begin{figure}[]
		\centering
\hspace{-3em}\includegraphics{./Figures/Legend_RS_split_out}\hspace{7em}
		\includegraphics{./Figures/Legend_Method_out}\\
  \hbox{
		\includegraphics{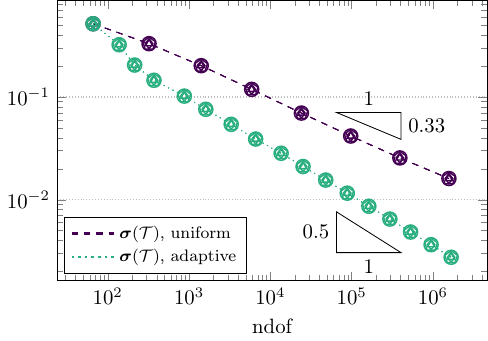}
		\includegraphics{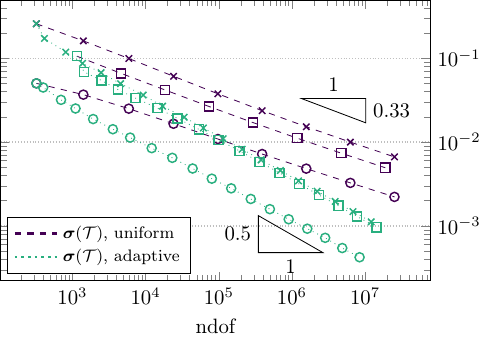}
		}
		\caption{Convergence history plot of the estimator $\boldsymbol{\sigma}(\T)$ for Morley FEM with different choices of $R, S\in\{id, J_h\}$ (left) and Morley,
			C0IP, dG I FEM (right)
		with unknown solution $u$ in Subsection \ref{sub:Numerical results for the Von Karman equation}}
		\label{fig:VKE_Dirac}
	\end{figure}

	\label{sub:Numerical results for the Von Karman equation}
	The second benchmark considers a point force $F=\delta_\zeta\in V^*$ located at the centroid $\zeta=(-1/6,-1/6)$ of the
	L-shaped domain.
 Example~\ref{ex:separation_condition} shows that the separation condition \eqref{eqn:single_load_assumption} holds for all newest-vertex refinements of the initial triangulation displayed in Figure \ref{fig:VKE_Dirac_meshes}.
 Given the localized estimator $\mu(\zeta,T)\coloneqq |T|^{1/2}$ if $T\in\T(\zeta)$ and $\mu(\zeta,T)\coloneqq0$ otherwise for the single load at $\zeta$,
 the discussion on single forces in Subsection \ref{sub:Single_forces} motivates
 the refinement indicator $\boldsymbol{\sigma}^2(T)$ for a triangle $T\in\T$ equal to
	\begin{align}\label{eqn:eta2}
		\mu&(\zeta, T) + |T|^2\Big(\big\|[u_h^{(1)},u_h^{(2)}]\big\|_{L^2(T)}^2 + \big\|[u_h^{(1)},u_h^{(1)}]\big\|_{L^2(T)}^2\Big) +\vartheta |T|^{1/2}\sum_{\ell=1,2}
					\big\|[\partial_{\nu\nu}u_h^{(\ell)}]\big\|_{L^2(\partial T\setminus\partial\Omega)}^2\\
     &+\sum_{\ell=1,2}\bigg(|T|^{1/2} \big\|[D^2u_h^{(\ell)}]_{\E}\tau_{\E}\big\|_{L^2(\partial T)}^2  +\sum_{E \in \mathcal{E}(T)} \Big(\big| \Pi_{E,0} \big[\partial_\nu u_h^{(\ell)}\big]_E \big|^2 + |T|^{-1}\hspace{-0.6em} \sum_{z \in \mathcal{V}(E)} \big|\big[u_h^{(\ell)}\big]_E(z)\big|^2 \Big) \bigg)\nonumber
	\end{align}
    with $\vartheta=1$ for $C^0$IP and $\vartheta=0$ else.
 	Figure \ref{fig:VKE_Dirac} displays optimal convergence rates of the adaptive algorithm driven by the refinement indicators \eqref{eqn:eta2}
	that improve on the observed suboptimal rate $1/3$ on uniformly refined meshes.
 Theorem \ref{thm:apost_vke_dirac} guarantees the observed equivalence of the unknown (undisplayed) error $\|{\bf e}_h\|_h\coloneqq\| {\bf u}-{\bf u}_h \|_{h}$ and the error estimator $\boldsymbol{\sigma}(\T)=\sqrt{\sum_{T\in\T}\boldsymbol{\sigma}^2(T)}$ up to vanishing oscillations.
	The convergence history plots for the different choices $R,S\in\{{\rm id}, I_\M,JI_\M\}$ overlap and are indistinguishable as
	highlighted for the Morley FEM on the left.
	Figure \ref{fig:VKE_Dirac_meshes} displays the adaptive refinement towards the atom of the point
	force $\zeta$ and an even stronger local refinement towards the reentrant corner.

	\subsection{Conclusive remarks}%
	\label{sub:Conclusion}

	The nested iteration of Newton's scheme for solving the nonlinear discrete problem is highly effective and reaches machine precision with 3 to 6 iterations in average.
 All variants with $R,S\in\{{\rm id},I_{\rm M},JI_{\rm M}\}$ lead to very similar accuracies.
 This is the first empirical confirmation of the overall equivalence of \cite[Thm.~8.3 and Thm.~9.1]{CCNNGRDS22}.
	While the theory requires a particular choice $S=J I_{\rm M}$ for the efficiency estimate in Theorem
	\ref{thm:releff}, 
	undisplayed computer experiments provide strong empirical evidence for efficiency of the presented a posteriori error estimator for \emph{any} choice of the
	operators $R,S \in \{{\rm id},I_{\rm M},JI_{\rm M}\}$.
	The classical schemes with $R=S={\rm id}$ are the easiest to implement
    and their a posteriori analysis may be combined with the reference scheme as suggested in Subsection~\ref{sub:Approximation from other discretisations} as a recommended overall strategy. The mandatory adaptive algorithm recovers optimal convergence rates in all examples and motivates future research on optimal convergence rates.

\section*{Acknowledgements}
The research of the first two authors has been supported by the Deutsche Forschungsgemeinschaft in the Priority Program 1748
under the project \emph{foundation and application of generalized mixed FEM towards nonlinear problems
in solid mechanics} (CA 151/22-2).  This paper has been supported by the  SPARC project 
(id 235) {\it the mathematics and computation of plates} and SERB POWER Fellowship SPF/2020/000019.
The second author is also supported by the \emph{Berlin Mathematical School, Germany}.
\bibliographystyle{amsplain}
{\footnotesize{\bibliography{NSEBib,references}}}
\newpage
\appendix
\renewcommand{\thesection}{\Alph{section}}
\renewcommand{\thesubsection}{\Alph{section}.\arabic{subsection}}
\title{{Supplement materials to the paper 'A posteriori error control for fourth-order semilinear problems with quadratic nonlinearity'}}
\author{Carsten Carstensen        
\quad\text{and}\quad Benedikt Gr\"a{\ss}le\quad\text{and}\quad
Neela Nataraj}
\date{}
\maketitle
This supplement contains three parts that provide further details on the practical application of the abstract results from Section \ref{sec:Abstract a posteriori error control} in Supplement A, an explicit construction of the test function $\psi$ used in the proof Theorem \ref{thm:apost_vke_dirac} for single forces in Supplement B, and a stopping criterion for solutions up to machine precision with Newton's method in Supplement C.

\section{Proof of Lemma~\ref{lem:uh-vh}}
The (finite) Taylor series expansion of $N_h$ at the root $u_h$ to $N_h(u_h)=0$ for the approximation of $N_h(v_h)$ provides
\begin{equation}\label{eqn:star}
DN_h(u_h; u_h-v_h)= \Gamma_{\pw} (R(u_h-v_h), R(u_h-v_h),S\,\bullet) -N_h(v_h) \in Y_h^*.
\end{equation}
By definition of the inf-sup constant \eqref{eqn:discreteinfsup}, there exists  $w_h \in Y_h$ with $\|w_h\|_{Y_h} \le 1$ 
and
\begin{align*}
\beta_h \|u_h-v_h\|_{X_h} &= DN_h(u_h; u_h-v_h,w_h) \nonumber \\
&= \Gamma_{\pw}(R(u_h-v_h), R(u_h-v_h),Sw_h)-N_h(v_h;w_h)
\nonumber \\
& \le \|\Gamma\|\|R\|^2\|S\|\| u_h -v_h\|^2_{X_h} + \|N_h(v_h)\|_{Y_h^*} \le \kappa \beta_h \| u_h -v_h\|_{X_h} + \|N_h(v_h)\|_{Y_h^*}
\end{align*}
with $\|\Gamma_\pw\|\|R\|^2\|S\|\|u_h-v_h\|_{X_h}\leq \kappa\beta_h$ in the last step.
This is the first assertion $(1-\kappa)\beta_h \|u_h-v_h\|_{X_h} \le  \|N_h(v_h)\|_{Y_h^*}$.
The second follows from the boundedness of    $DN_h(u_h)$ and \eqref{eqn:star}; in fact
\begin{align*}
\|N_h(v_h)\|_{Y_h^*} & = \| \Gamma_{\pw}(R(u_h-v_h), R(u_h-v_h),S\,\bullet) - DN_h(u_h; u_h-v_h)\|_{Y_h^*} \\
& 
\le \left( \| \Gamma\|\|R\|^2\|S\| \|u_h-v_h\|_{X_h} +   \|DN_h(u_h)\|_{X_h^* \times Y_h^*}\right) \| u_h -v_h\|_{X_h} \\
& \le \left(  \kappa \beta_h  +   \|DN_h(u_h)\|_{X_h^* \times Y_h^*}\right)\| u_h -v_h\|_{X_h} . \qquad \qquad\qquad \qquad \qquad \qquad \qquad \qquad \hfill \qed
\end{align*}
\section{Design of the test function $\psi$}
Given a point $\zeta\in \mathrm{int}(E)$ on the interior of an edge $E$ of the triangulation $\T$, the construction a test function $\psi$ with $\psi(\zeta)=1$ and a list of orthogonalities \eqref{eqn:psi_def} used in the proof of Theorem~\ref{thm:apost_vke_dirac} based on one-dimensional Jacobi polynomials follows in three steps.

\medskip
\noindent{\it Step 1} discusses the orthogonal Jacobi polynomials $P_n^{(4,4)} \in P_n[-1,1]$ of degree $n \in {\mathbb N}_0$ that reflect the weight $\rho(x):=(1-x^2)^4$ for $-1 \le x \le 1$. 
The well-known three-term recurrence relation reveals
\[c_n\coloneqq P_{2n}^{(4,4)}(0) = (-4)^{-n}\begin{pmatrix}
    2n + 4\\n
\end{pmatrix}\qquad\text{for }n\in\mathbb N_0\]
and guarantees $P_{2n}(0) \neq 0$. Given $k \in {\mathbb N}_0$, select $n\in {\mathbb N}$ with $2n \ge k$ and rescale to define $\psi\coloneqq P_{2n}^{(4,4)}/c_n\in P_{2n}[-1,1]$ with $\psi(0)=1$ and 
$\rho \psi \perp P_k[-1,1]$ in $L^2(-1,1)$. 
Observe that $w(x,y)=(1-x^2)^2 (1-y^2)^2$ implies $w(x,x) =\rho(x)$.
The polynomial $f\in P_{8+2n}(Q)$ defined by
\[f(x,y):= w(x,y) \psi({(x+y)}/{2})\qquad\text{for } x,y \in[-1,1] \] 
 on the cube $Q=[-1,1]^2$ satisfies along the diagonal $D:=\{(x,x): -1 \le x \le 1 \}$ that
\[f(0,0)=1 \text{ and } fq \perp P_n(D) \text{ in } L^2(D)^2 \text{ for all }q_k \in P_k(D)^2.\]
 By symmetry of $w$ along $D$, $\nabla w|_D \perp \nu_D=(1;-1)/\sqrt{2}$ pointwise along $D$. Since the gradient $\psi'({(x+y)}/{2})(1/2;1/2)$ of $(x,y) \mapsto \psi({(x+y)}/{2})$ is also perpendicular to $\nu_D$ along $D$, we infer $\nu_D \cdot \nabla f|_D=0$. The scaling by $h>0$ leads to $g(x,y)=f(x/h,y/h)$ with
\begin{align*}
& g \in P_{8+2n}(Q) \cap H^2_0(Q), \; {\rm supp} \; g \subseteq Q, \; g \perp P_k(hD) \; {\rm  in} \; L^2(hD),\\
& \nabla g \cdot \nu_D=0 \text{ on } hD, \text{ and }|g|_{H^s(hQ)}=h^{1-s} |f|_{H^s(Q)} \text{ for } s=0,1,2.
\end{align*}
Those properties are inherited by transformations in another Cartesian coordinate system (by translation and rotation).

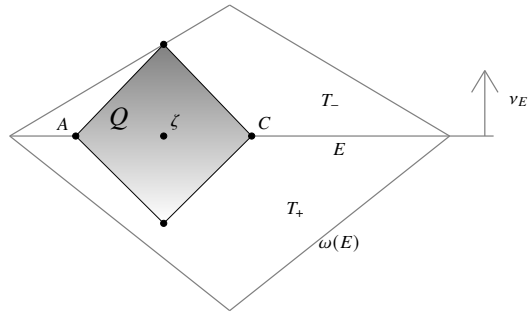
\begin{figure}[h]
\begin{center}
\resizebox{0.4\textwidth}{!}{
\begin{tikzpicture}
\draw[gray, thick] (-5,0) -- (0,3);
\draw[gray, thick] (0,3) -- (5,0);
\draw[gray, thick] (5,0) -- (0,-4);
\draw[gray, thick] (0,-4) -- (-5,0);
\draw[gray, thick] (-5,0) -- (6,0);
\shadedraw (-3.5,0) -- (-1.5,2.1) -- (0.5,0) -- (-1.5,-2) -- (-3.5,0);
\filldraw[black] (-3.5,0) circle (2pt) node[anchor=south east]{$A$};
\filldraw[black] (0.5,0) circle (2pt) node[anchor=south west]{$C$};
\filldraw[black] (-1.5,0) circle (2pt) node[anchor=south west]{$\zeta$};
\filldraw[black] (-1.5,-2) circle (2pt);
\filldraw[black] (-2.5,0) node [anchor=south]{\Large{$Q$}};
\filldraw[black] (-1.5,2.1) circle (2pt);
\node at (2.5,0) [anchor=north]{$E$};
\node at (2.3,1.1) [anchor=north]{$T_-$};
\node at (1.5,-2) [anchor=south]{$T_+$};
\node at (2.5,-2.8) [anchor=south]{$\omega(E)$};
\draw[gray, thick] (5.8,0) -- (5.8,1.5);
\draw[gray, thick] (5.5,1) -- (5.8,1.5);
\draw[gray, thick] (5.8,1.5) -- (6.1,1);
\node at (6.2,0.8) [anchor=west]{$\nu_E$};
\end{tikzpicture}
}
\end{center}
\caption{Triangles $T_+$, $T_-$, $\omega(E)$, and square $Q \subset T_+ \cup T_-$} \label{suppl:fig}
\end{figure}

\medskip \noindent \noindent{\it Step 2} constructs an edge bubble-function over the edge-patch.  Given an interior edge $E= \partial T_+ \cap \partial T_-$ shared by the triangles $T_\pm \in \T$ and patch $\omega(E)$ and $\zeta \in E$ with ${\rm dist} (\zeta,v) \approx h_E \approx h_{T_+} \approx h_{T_-}$ (from shape-regularity of $\T$), let $Q$ be the maximal square with edge-size $h>0$ and midpoint $\zeta$ that belongs to $T_+ \cup T_-  \supset Q$ such that one diagonal ${\rm conv} \{A,C\} \subset E$ lies on the edge $E$  as displayed in Figure~\ref{suppl:fig}.
\begin{comment}
\begin{figure}[h]
\begin{center}
\resizebox{0.4\textwidth}{!}{
\begin{tikzpicture}
%
\draw[gray, thick] (-5,0) -- (0,3);
\draw[gray, thick] (0,3) -- (5,0);
\draw[gray, thick] (5,0) -- (0,-4);
\draw[gray, thick] (0,-4) -- (-5,0);
%
\draw[gray, thick] (-5,0) -- (6,0);
%
\draw[gray, thick] (-3,0) -- (-1.5,2.1);
\draw[gray, thick] (-1.5,2.1) -- (0,0);
\draw[gray, thick] (0,0) -- (-1.5,-2);
\draw[gray, thick] (-1.5,-2) -- (-3,0);
%
\filldraw[black] (-3,0) circle (2pt) node[anchor=south east]{A};
\filldraw[black] (0,0) circle (2pt) node[anchor=south west]{C};
\filldraw[black] (-1.5,0) circle (2pt) node[anchor=south west]{$\zeta$};
\filldraw[black] (-1.5,-1)  node [anchor=south]{Q};
\filldraw[black] (-1.5,2.1) circle (2pt);
%
\node at (2.5,0) [anchor=north]{$E$};
\node at (2.3,1.1) [anchor=north]{$T_-$};
\node at (1.5,-2) [anchor=south]{$T_+$};
\node at (2.5,-2.8) [anchor=south]{$\omega(E)$};
%
\draw[gray, thick] (5.8,0) -- (5.8,1.5);
%
\draw[gray, thick] (5.5,1) -- (5.8,1.5);
\draw[gray, thick] (5.8,1.5) -- (6.1,1);
%
\node at (6.2,0.8) [anchor=west]{$v_E$};
\end{tikzpicture}}
\end{center}
\caption{Triangles $T_+$, $T_-$, $\omega(E)$, and $Q$} \label{suppl:fig}
\end{figure}
\end{comment}
\medskip \noindent A translation by $\zeta$ and a rotation to fit $\zeta, A, C \in E$ leads to a function $g \in H^2_0(Q) \cap P_{8+2n}$ as designed in Step 1 with scaling $|g|_{H^s(Q)} \approx h^{1-s}$ for $s=0,1,2$ and with $g(\zeta) =1$ and various orthogonalities.

\medskip \noindent \noindent{\it Step 3} is the final design of $\psi \perp P_k(\T)$. This is more standard than the previous design steps with a cubic bubble-function $b_{T_+}$ and  
$b_{T_-}$ and polynomials $q_{\pm} \in P_k(T_{\pm})$ such that $g- b^2_{T_{\pm}} q_{\pm} \perp P_k(T_{\pm})$ in $L^2(T_{\pm})$. This leads to the function $\psi:=g-b^2_{T_+}q_+ -b^2_{T_-}q_- \in H^2_0(Q)\subset V $ with all the desired orthogonality conditions, ${\rm supp} \; \psi \subset T_+ \cup T_-$, and $\psi(\zeta)=1$. The scaling $\trinl \psi \trinr \lesssim h^{-1} \approx h_E^{-1}$ follows from that of $g$ and the following routine estimate for $q_{\pm}$. Inverse estimates show 
\[ \|q_{\pm} \|_{L^2(T_{\pm})}^2 \approx \|b_{T_{\pm}}q_{\pm} \|_{L^2(T_{\pm})}^2 =\int_{T_{\pm}} gq \dx \le  \|g\|_{L^2(T_{\pm})} \|q_{\pm}\|_{L^2(T_{\pm})}. \] 
This, two Friedrichs's inequalities, and an inverse inequality conclude the proof with
\[ \|b_{T_{\pm}}q_{\pm} \|_{L^2(T_{\pm})} \approx
\|q_{\pm}\|_{L^2(T_{\pm})} \le
\|g_{\pm}\|_{L^2(T_{\pm})}
\le{\pi^{-2}} {h_{T_{\pm}}} \trinl g \trinr.\qed\]

\section{Accurate solution to semilinear problems}
This supplement provides algorithmic details on the adaptive computations in Section \ref{sec:Numerical experiments} with particular focus on the implementation of Newton's method for the solution of the discrete equation \eqref{eqn:discrete} up to machine precision controlled by the termination criterion of Lemma~\ref{lem:uh-vh}.

	\subsection*{Accurate approximations with Newton's method}%
	\label{ssub:Approximation of the discrete solution}
	Starting from the initial guess $u_{\ell}^0\coloneqq I_hJ I_\M u_{\ell-1}$ for nested iteration with initialisation $u_0^0\coloneqq0$ on the coarsest mesh, our implementation computes the Newton
	iterates $u_{\ell}^{k+1}\coloneqq u_{\ell}^{k}-DN_h(u_\ell^k)^{-1}N_h(u_\ell^k)$ with the LU decomposition in Julia for an exact Newton update.
    In the present situation of Section \ref{sec:Numerical experiments}, $a_h$ is a scalar product associated to the linear operator $A_h\in L(V_h;V_h)$ and induces the method-dependent norm $\|\bullet\|_a\coloneqq a_h(\bullet,\bullet)^{1/2}$ in $V_h$. 
    Let $u_\ell$ denote the exact discrete solution and recall the equivalence of the algebraic error to the residual from Lemma~\ref{lem:uh-vh}. This and the Riesz isomorphism  
    \begin{align*}
        \big\|A_h^{-1}N_h\big(u_\ell^{k+1}\big)\big\|_a= \left\|N_h\big(u_\ell^{k+1}\big)\right\|_{a,*}\coloneqq\sup_{v_h\in
		V_h}\frac{\left|N_h\left(u_\ell^{k+1}\right)v_h\right|}{\|v_h\|_a}\approx \big\|u_\ell-u_\ell^{k+1}\big\|_a
    \end{align*}
    motivates a stopping criterion on the computable quantity $\big\|A_h^{-1}N_h\big(u_\ell^{k+1}\big)\big\|_a$ in two stages.
	The first step iterates until
	\begin{align}\label{eqn:stopp_crit}
		\big\|A_h^{-1}N_h\big(u_\ell^{k+1}\big)\big\|_a\leq
		tol\;\left(\left\|u_\ell^{k+1}\right\|_a+\left\|u_\ell^k\right\|_a\right)
	\end{align}
	holds with $tol\coloneqq10^{-4}$ in the benchmarks of Section \ref{sec:Numerical experiments}.
Once this coarse condition is satisfied
	the algebraic error is considered sufficiently small such that quadratic
	convergence can be expected through the
	Newton-Kantorovich theorem.
	The second stage computes further iterates until
	\begin{align}\label{eqn:stopp_crit2}
		\big\|A_h^{-1}N_h\big(u_\ell^{k}\big)\big\|_a\leq\big\|A_h^{-1}N_h\big(u_\ell^{k+1}\big)\big\|_a.
	\end{align}
    This suggests the approximate solution $u_\ell= u_\ell^k$ is accurate up to machine precision. In average, the benchmarks in Section \ref{sec:Numerical experiments} (with nested iteration) perform 1 to 2 iterations until \eqref{eqn:stopp_crit} and another 0 to 5 iterations until \eqref{eqn:stopp_crit2} holds.
    \begin{rem}[General bilinear forms $a_h$]
    If $a_h$ does not define a scalar product in $V_h$,
    Example~\ref{ex:computation_residual} and Lemma~\ref{lem:uh-vh} provide an alternative approach for the computation of the algebraic residual $$\big\|A_h^{-1}N_h\big(u_\ell^{k+1}\big)\big\|_h\approx \left\|N_h\big(u_\ell^{k+1}\big)\right\|_{*}\approx\big\|u_\ell-u_\ell^{k+1}\big\|_h$$ in the norm $\|\bullet\|_h$ from \eqref{hnorm}.
    \end{rem}

\end{document}